\documentclass[11pt]{amsart}

\usepackage[T1]{fontenc}
\usepackage[utf8]{inputenc}

\usepackage{hyperref}

\allowdisplaybreaks[4]
\usepackage{algorithm}
\usepackage{algpseudocode}

\usepackage{amsmath,mathtools}
\usepackage{amssymb,bbm}
\usepackage{mathrsfs}
\usepackage{amsthm}
\usepackage{bm}
\usepackage{graphicx}
\usepackage{tikz}
\usepackage{epstopdf}
\usepackage{booktabs}
\usepackage{epsfig}

\usepackage[caption=false]{subfig}
\usepackage{caption}

\usepackage{float}

\newfloat{longfigure}{htbp}{lof}
\floatname{longfigure}{Figure}
\newsubfloat{longfigure}

\usepackage{geometry}
\usepackage{color}

\usepackage{color}
\theoremstyle{plain}
\newtheorem{lemma}{Lemma}[section]
\newtheorem{theorem}[lemma]{Theorem}

\newtheorem{remark}[lemma]{Remark}

 \makeatletter
\@namedef{subjclassname@2020}{%
  \textup{2020} Mathematics Subject Classification}
\makeatother
\numberwithin{equation}{section}

\begin{document}
\title[]{Weak error estimates of a spectral Galerkin approximation to the stochastic Burgers equation driven by rough white noise}
\author{Minxing Zhang, Yongkui Zou}
\address{School of Mathematics, Jilin University, Changchun 130012, China}
\email{zhangminxingjlu@163.com, zouyk@jlu.edu.cn}

\author{Yanzhao Cao}
\thanks{Yanzhao Cao is the corresponding author.}
\address{Department of Mathematics and Statistics, Auburn University, Auburn, AL 36840, USA}
\email{yzc0009@auburn.edu}
\thanks{This work is supported by National Natural Science Foundation of China (12571428) and the Project for Enhancing the Research and Innovation Capacity of the Education Department of Jilin Province (2025KC028).}
\begin{abstract}
 	In this paper, we investigate weak error estimates for a spectral Galerkin approximation to a stochastic Burgers equation driven by rough space-time white noise.
 	The roughness of the driving noise results that the mild solution possesses only $L^\infty$-regularity rather than any positive-order differentiability order regularity. 
 	The theory of classical numerical analysis reveals that the convergence rate of numerical solution is determined by the smoothness of the exact solution, and hence $L^\infty$-regularity is insufficient to yield any positive convergence rate. 
 	Indeed, De Vecchi et al.\ \cite{DeVecchi2026} studied a spectral Galerkin approximation to such equation and only proved pathwise convergence.
 	However, we aim to derive positive weak convergence rate for the spectral Galerkin approximation, which is determined by the roughness order of noise.
 	To overcome the difficulty caused by the absence of smoothness of the exact solution, we turn to study the differentiability regularity of the first and second order Fr\'echet derivatives of the solution  to a Kolmogorov equation associated with the spectral Galerkin approximation. We manage to prove that they possess nearly second order smoothness, enabling us to derive a positive weak convergence rate for the spectral Galerkin approximation.	Finally, we conduct a numerical experiment to illustrate the theortical analysis.
\end{abstract}
\keywords {Stochastic Burgers equation, rough space-time white noise, Kolmogorov equation, weak error estimates, spectral Galerkin method}
\subjclass[2020]{65C30, 60H15, 60H35}

\maketitle

\section{Introduction}
\label{Sect1}

In recent years, numerical approximations to  nonlinear stochastic partial differential equations (SPDEs) with non-globally Lipschitz continuity, such as stochastic Cahn–Hilliard equations \cite{CuiJianbo2020,CuiJianbo2021_2,CaiMeng2023,HongJialin2024}, stochastic Allen–Cahn equations \cite{ChenChuchu2025,ZhangZou2026JSC,QiXiao2026JSC,ZhangZou2027JCAM}, and stochastic Burgers equations \cite{Jentzen2013,WangYibo2026} have been extensively studied to explore the statistical properties of the underlying mathematical models, where both strong and weak error estimates have been used to investigate the numerical convergence. The strong error characterizes the average deviation of a numerical trajectory from the exact one \cite{LiuZhihui2021,CuiJianbo2021_2}, while the weak error measures the discrepancy in sampling the probability law of the exact solution \cite{CuiJianbo2019,Brehier2026SCHE}.

The Burgers equation is a fundamental model in fluid dynamics, describing the interaction between nonlinear convective effects and viscous dissipation.
A stochastic Burgers equation (SBEs) incorporates random forcing and serves as a simplified model for studying turbulence and the interaction between nonlinear dynamics and stochastic fluctuations \cite{Gyongy2006}.
Numerical approximations to SBEs have been extensively studied in the literature.
Early numerical analyses were primarily concerned with pathwise approximation  and strong convergence.
Alabert and Gy\"ongy~\cite{Gyongy2006} constructed a finite difference spatial
approximation to SBEs driven by additive space-time white noise and established almost sure convergence orders. Bl\"omker and Jentzen~\cite{Jentzen2013} applied a Galerkin approximation to SBE driven by additive white noise and proved a pathwise convergence rate arbitrarily close to $1/2$ in the maximum norm.
Bl\"omker et al.\ \cite{Blomker2013} subsequently applied a fully discrete exponential scheme to SBEs driven by correlated additive noise and obtained pathwise convergence rates. 
Hairer et al.\ \cite{Hairer2014} applied a Fourier multipliers combining finite difference scheme to establish a spatial semi-discrete approximation to SBEs driven by multiplicative space-time white noise and investigated pathwise convergence rate of order $1/6$. Hairer and Matetski~\cite{Hairer2016} later improved the pathwise convergence rate to almost $1/2$ in additive noise case.
Jentzen et al.\ \cite{Jentzen2020_2} applied a spectral Galerkin method and nonlinearity-truncated exponential scheme to SBEs driven by additive space-time white noise and derived strong convergence properties.
Hutzenthaler and Jentzen~\cite{Hutzenthaler2020} applied a spectral Galerkin method to SBEs driven by multiplicative noise and proved strong convergence. 
Hutzenthaler and Link~\cite{Hutzenthaler2022} applied a spectral Galerkin method and tamed exponential Euler scheme to SBEs driven by multiplicative trace-class noise and proved that the spatial and temporal strong convergence orders are of $1$ and $1/2$, respectively.
Wang and Cao~\cite{WangYibo2025} applied a spectral Galerkin method and tamed exponential Euler scheme to SBEs driven by additive fractional noise and obtained strong convergence properties.

More recently, there have been a few studies focusing on the weak convergence analysis of numerical approximations to SBEs \cite{Brehier2026SBE,ZhangZou2026DCDSb,LiBinjie2026}. Br\'ehier et al.\ \cite{Brehier2026SBE} applied a spectral Galerkin method to SBEs driven by additive trace-class noise, and proved nearly the second order weak convergence.  Their proof ultilizes a special property that the time integral of the squared $H^1$-norm of the spectral Galerkin mild solution is exponential integrable.
They also pointed out in the same paper that the weak error estimates for SBEs driven by space-time white noise "is an open problem". 

In this paper, we investigate weak error estimates for a spectral Galerkin approximation to SBEs driven by rough space-time white noise $(-\Delta)^{\gamma}dW(t)$ with roughness order $\gamma\in[0,\frac{1}{4})$. 
By applying the Da Prato-Debussche trick \cite{DaPrato2007}, we only obtain $L^\infty$ a priori estimates for the exact mild solution. The classical numerical analysis indicates that the strong convergence order is at most equal to the smoothness order of the exact mild solution, whereas the $L^\infty$-integrability is impossible to provide any positive convergence order. Indeed, De Vecchi et al.\ \cite{DeVecchi2026} studied a spectral Galerkin approximation to such equation and only proved pathwise convergence property without getting any positive convergence order. 
Nevertheless, we still aim to derive a positive  weak convergence order for the spectral Galerkin approximation in this work. The main contribution is that we find the regularity of the  solution to a Kolmogorov equation associated with the spectral Galerkin approximation can provide extra weak convergence order. Therefore,  we manage to establish a priori estimates for the first- and second-order Fr\'echet derivatives of the solution to Kolmogorov equation, as well as their composition with nearly $1$ power of the Laplacian. Finally, we  prove that the weak convergence rate  is arbitrarily close to $1-4\gamma$.

The rest of this paper is organised as follow. In Section \ref{Sect2}, we introduce SBEs and derive the $L^\infty$-estimates of its solution. We also construct a spectral Galerkin approximation to SBEs and collect some regularities of its solution. In Section \ref{Sect3}, we derive a Kolmogorov equation associated with the spatial semi-discrete scheme and study some Sobolev regularity estimates of its solution. In Section \ref{Sect4}, we investigate weak error estimates. In Section \ref{Sect5}, we conduct a numerical experiment to illustrate the theortical analysis.

\section{Preliminaries}
\label{Sect2}

\subsection{Working spaces}

Let $(U,(\cdot, \cdot)_{U})$ and $(V,(\cdot, \cdot)_{V})$ be two real separable Hilbert spaces. By $\mathcal{L}(U;V)$ we denote the Banach space of all bounded linear operators from $U$ to $V$ equipped with the operator norm. For simplicity, we write $\mathcal{L}(U) = \mathcal{L}(U;U)$.

For any $T\in\mathcal{L}(U;V)$, denote by $\{\sigma_k\}_{k\in\mathbb{N}}$ the eigenvalues of $(TT^*)^{\frac{1}{2}}$. For any $1\le p<\infty$, define a Banach space of linear operators
\begin{equation*}
	\mathcal{L}_p(U;V) = \big\{T\in\mathcal{L}(U;V):\ \|T\|_{\mathcal{L}_p(U;V)} = \big(\sum_{k=1}^{\infty}\sigma_k^p\big)^{\frac{1}{p}}<\infty\big\}.
\end{equation*}
Especially, $\mathcal{L}_1(U;V)$ and $\mathcal{L}_2(U;V)$ are the spaces of trace-class operators and Hilbert-Schmidt operators, respectively. To simplify notation, we write $\mathcal{L}_p(U)=\mathcal{L}_p(U;U)$. 
Let $\{e_k\}_{k=1}^{\infty}$ be an orthonormal basis in $U$, then we have for $T\in \mathcal{L}_1(U)$,
\begin{equation*}
	\mathrm{Tr}(T) = \sum_{k=1}^{\infty}(Te_{k}, e_{k})_{U} \le \|T\|_{\mathcal{L}_{1}(U)}.
\end{equation*}
Furthermore, for any $S\in \mathcal{L}(U;V)$, $T_1,T_2 \in \mathcal{L}_2(U)$, there hold \cite{DaPrato2014,Kruse2014}
\begin{equation}\label{eq7.2.3}
	\begin{aligned}
		&\|ST_1\|_{\mathcal{L}(U;V)}\le\|ST_1\|_{\mathcal{L}_2(U;V)}\le\|S\|_{\mathcal{L}(U;V)}\|T_1\|_{\mathcal{L}_2(U)},\\ &\|T_1T_2\|_{\mathcal{L}_2(U)}\le\|T_1T_2\|_{\mathcal{L}_1(U)}\le\|T_1\|_{\mathcal{L}_2(U)}\|T_2\|_{\mathcal{L}_2(U)}.
	\end{aligned}
\end{equation}

For any integer $m\ge 1$, by $\mathcal{C}_b^m(U)$ denote the space of not necessarily bounded functions $\phi: U \rightarrow \mathbb{R}$ having continuous and bounded Fr\'echet derivatives $D^i\phi(u)$ $(1\le i\le m)$ equipped with a maximum seminorm $|\phi|_{\mathcal{C}_b^m(U)}$.
Noticing that the first order Fr\'echet derivative $D\phi(u)$ of $\phi$ at $u\in U$ is a bounded linear functional on $U$ and according to Riesz representation theorem, there exists a unique element in $U$, still denoted by $D\phi(u)$, such that
\begin{equation}\label{eq7.2.4}
	D\phi(u)h = ( D\phi(u),h)_U,\quad |D\phi(u)h| \le |\phi|_{\mathcal{C}_b^1(U)}\|h\|_U, \quad \forall u,h \in U.
\end{equation}

Let $\mathcal{D}=(0,1)$. Denote by $W^{k,p}$ $(k\ge 0,p\ge 1)$ the standard Sobolev spaces defined on $\mathcal{D}$, and simply write $H^k=W^{k,2}$. Especially, $H = W^{0,2}(\mathcal{D}) = L^2(\mathcal{D})$. The Sobolev embedding $W^{s,p}\hookrightarrow W^{r,q}$ holds for $s\ge r$ and $s-\frac{1}{p}\ge r-\frac{1}{q}$. Let $V = L^{\infty}(\mathcal{D})$ be the space of essentially bounded functions with norm $\|\cdot\|_V$. 

Let $A =-\Delta:\mathrm{Dom}(A) = H_0^1 \cap H^2 \subset H \to H$ be a Laplacian operator with homogeneous Dirichlet boundary condition. 
Then, the operator $A$ possesses a family of eigenvalues $\lambda_i = i^2\pi^2$ with complete orthonormal eigenbases $\varphi_i(x)=\sqrt{2}\sin(i\pi x)$. 
For any $\alpha \ge 0$, define a fractional power of $A$ by $\displaystyle A^{\alpha}v = \sum_{j=1}^{\infty}\lambda_j^\alpha(v,\varphi_j)\varphi_j$ and a space $\displaystyle\dot{H}^{\alpha} = \{v \in H: \| v \|_{\alpha} =\| A^{\frac{\alpha}{2}}v\| < \infty\}$. Then $(\dot{H}^{\alpha},\| \cdot \|_{\alpha})$ is a Hilbert space. The operator $-A$ generates an analytic semigroup $S(t) = e^{-tA}$ on $H$. It is pointed out in \cite{Thomee2006} that there exists a constant $C>0$ such that
\begin{equation}\label{eq3.2.2}
	\begin{aligned}
		A^{\alpha}S(t) = S(t)A^{\alpha} \ \text{ on } \dot{H}^{2\alpha}  \text{ and } 	\|A^{\alpha}S(t)\|_{\mathcal{L}(H)} \le Ct^{-\alpha}, \quad \forall 0<t\le T,\ \alpha\ge 0,\\
		\int_s^t\|A^{\alpha}S(t-r)u\|^2dr \le C(t-s)^{1-2\alpha}\|u\|^2,\quad \forall u\in H,\ 0\le s\le t\le T,\ 0\le\alpha\le\frac{1}{2}.
	\end{aligned}
\end{equation}
Furthermore, the semigroup $S(t)$ satisfies some contractivity properties: for any $1\le p\le q<\infty$ and $u\in L^p(\mathcal{D})$,
\begin{equation}\label{eq11.2.3}
	\begin{aligned}
		\|S(t)u\|_{L^q(\mathcal{D})} \le&\ Ct^{-\frac{1}{2}(\frac{1}{p}-\frac{1}{q})}\|u\|_{L^p(\mathcal{D})},\\
		\|S(t)u\|_{V} \le&\ Ct^{-\frac{1}{2p}}\|u\|_{L^p(\mathcal{D})}.
	\end{aligned}
\end{equation}
Throughout this paper, by $\partial_x u$ we denote  the weak derivative of $u$ with respect to the spatial variable $x$. According to \cite{Brehier2026SBE}, for any $p>1$ and $\delta>0$, there exists a constant $C>0$ such that
\begin{equation}\label{eq12.2.3}
	\begin{aligned}
		\|A^{-\frac{1}{2}}\partial_x u\|_{L^p(\mathcal{D})} \le&\ C\|u\|_{L^p(\mathcal{D})},\quad \forall u\in W^{1,p}(\mathcal{D}),\\
		\|A^{-\frac{1}{4}-\delta}u\| \le&\ C\|u\|_{L^1(\mathcal{D})},\quad\forall u\in L^1(\mathcal{D}).
	\end{aligned}
\end{equation}

\subsection{Stochastic Burgers equation and properties}

Let $(\Omega,\mathcal{F},\mathbb{P})$ be a completed probability space. Fix $T>0$ and consider a stochastic Burgers equation driven by rough noise
\begin{equation}\label{eq14.1.1}
	\begin{aligned}
		dX(t)+A X(t)dt =&\ F(X(t))dt+A^{\gamma}dW(t),\ t\in (0,T],\\
		\quad X(0) =&\ \xi,
	\end{aligned}
\end{equation}
where $F: L^4(\mathcal{D}) \rightarrow \dot{H}^{-1}$ is a nonlinear mapping defined by
\begin{equation*}
	{}_{\dot H^{-1}}\langle F(u),v\rangle_{\dot H^1} :=
	-\int_{\mathcal{D}}u^2(x)\partial_xv(x)dx, \quad u\in L^4(\mathcal{D}),\ v\in\dot H^1.
\end{equation*}
Equivalently, $F(u)=\partial_x(u^2)$ in the sense of distributions.
If $u$ is sufficiently regular, this coincides with $F(u)=\nabla (u^2) = 2u\nabla u$. Moreover, $W(t)$ is a space-time white noise.
By $X(t,\xi)$ we denote the solution of \eqref{eq14.1.1} satisfying $X(0,\xi) = \xi$.
We assume

\begin{itemize}
	\item[(A1)] $\gamma\in[0,\frac{1}{4})$.
	\item[(A2)] The initial value $\xi\in \dot{H}^1$ is a deterministic function.
\end{itemize} 

We emphasize that if $\gamma\ge\frac{1}{4}$, the stochastic convolution exists only in the sense of distributions \cite{Hairer2013}, and such a case is not considered in this work. If $\gamma<-\frac{1}{2}$, $A^{\gamma}$ becomes a trace-class operator and this case is already studied by Br\'ehier et al. \cite{Brehier2026SBE}.

In next lemma, we study some a priori estimates of $F$.

\begin{lemma}\label{l11.2.1}
	The nonlinear mapping $F$ satisfies
	\begin{equation*}
		(u,F(u)) = 0,\quad \forall u \in \dot{H}^{1}.
	\end{equation*}
	Moreover, for any $v,w\in\dot{H}^1$ and $p>1$, there exists a constant $C>0$ such that
	\begin{align*}
		(v,DF(u)v) \le&\ \|u\|_1\|v\|^2_{L^4(\mathcal{D})},\\
		\|A^{-\frac{1}{2}}DF(u)v\| \le&\ 2\|u\|_V\|v\|,\\ 
		\|A^{-\frac{1}{2}}D^2F(u)(v,w)\|_{L^1(\mathcal{D})} \le&\ C\|w\|_{L^{2p}(\mathcal{D})}\|v\|_{L^{2p}(\mathcal{D})}.
	\end{align*}
\end{lemma}

\begin{proof}
	From the integration by parts formula and the homogeneous Dirichlet boundary condition, it follows that for any $u\in\dot{H}^1$,
	\begin{equation*}
		(u,F(u))=\int_{\mathcal{D}}u(x)\partial_x(u^2(x))dx=-\int_{\mathcal{D}}\partial_xu(x)u^2(x)dx.
	\end{equation*}
	On the other hand, we have
	\begin{equation*}
		(u,F(u))=2\int_{\mathcal{D}}u^2(x)\partial_xu(x)dx=-2(u,F(u)).
	\end{equation*}
	Thus we deduce $(u,F(u))=0$.
	
	Since $\dot{H}^1$ is a Banach algebra, there exists a constant $C>0$ such that $\|\partial_x(v^2)\|=\|v^2\|_1\le C\|v\|_1^2$. Therefore, for $u,v\in\dot{H}^1$,
	\begin{align*}
		&\ \|F(u+v)-F(u)-2\partial_x(uv)\|=\|\partial_x[(u+v)^2-u^2]-2\partial_x(uv)\|\\
		=&\ \|\partial_x(v^2)\|\le C\|v\|_1^2.
	\end{align*}
	Moreover, $\|2\partial_x(uv)\|\le C\|u\|_1\|v\|_1$. Consequently,
	\begin{equation*}
		\lim_{\|v\|_1\to0}\frac{\|F(u+v)-F(u)-2\partial_x(uv)\|}{\|v\|_1}=0,
	\end{equation*}
	which implies that the restriction $F:\dot{H}^1\to H$ is Fr\'echet differentiable and
	\begin{equation}\label{eq14.2.5}
		DF(u)v=2\partial_x(uv)=2\partial_xu\cdot v+2u\partial_xv.
	\end{equation}
	Repeatedly applying the integration by parts formula, we obtain
	\begin{align*}
		(v,DF(u)v)=&\ 2\int_{\mathcal{D}}v(x)\partial_x(u(x)v(x))dx=-2\int_{\mathcal{D}}\partial_xv(x)u(x)v(x)dx\\
		=&\ -\int_{\mathcal{D}}u(x)\partial_x(v^2(x))dx=\int_{\mathcal{D}}\partial_xu(x)v^2(x)dx\\
		\le&\ \|\partial_xu\|\|v^2\|=\|u\|_1\|v\|_{L^4(\mathcal{D})}^2.
	\end{align*}
	
	By the dual characterization of the $\dot{H}^{-1}$-norm, the integration by parts formula, and the H\"older inequality, we have
	\begin{align*}
		&\ \|A^{-\frac{1}{2}}DF(u)v\|= \|DF(u)v\|_{-1}=\sup_{\substack{z\in\dot{H}^1\\ \|z\|_1=1}}|(DF(u)v,z)|\\
		=&\ 2\sup_{\substack{z\in\dot{H}^1\\ \|z\|_1=1}}\Bigl|\int_{\mathcal{D}}u(x)v(x)\partial_xz(x)dx\Bigr|
		\le 2\sup_{\substack{z\in\dot{H}^1 \|z\|_1=1}}\|u\|_V\|v\|\|z\|_1=2\|u\|_V\|v\|.
	\end{align*}
	
	Notice that for any $u,v,w\in\dot{H}^1$,
	\begin{equation*}
		DF(u+w)v-DF(u)v=2\partial_x(wv),
	\end{equation*}
	and $\|2\partial_x(wv)\|\le C\|w\|_1\|v\|_1$. Hence
	\begin{equation*}
		\sup_{\substack{v\in\dot{H}^1\\ \|v\|_1=1}}\|DF(u+w)v-DF(u)v-2\partial_x(wv)\|=0.
	\end{equation*}
	It follows that $F:\dot{H}^1\to H$ is twice continuously Fr\'echet differentiable and
	\begin{equation*}
		D^2F(u)(v,w)=2\partial_x(vw).
	\end{equation*}
	By \eqref{eq12.2.3}, the Sobolev embedding $L^p(\mathcal{D})\hookrightarrow L^1(\mathcal{D})$, and the H\"older inequality, we obtain for any $p>1$,
	\begin{align*}
		&\ \|A^{-\frac{1}{2}}D^2F(u)(v,w)\|_{L^1(\mathcal{D})}= 2\|A^{-\frac{1}{2}}\partial_x(vw)\|_{L^1(\mathcal{D})} \le 2\|A^{-\frac{1}{2}}\partial_x(vw)\|_{L^p(\mathcal{D})}\\
		\le&\ C(p)\|vw\|_{L^p(\mathcal{D})}\le C(p)\|v\|_{L^{2p}(\mathcal{D})}\|w\|_{L^{2p}(\mathcal{D})}.
	\end{align*}
	Thus the proof is completed.
\end{proof}

The following lemma introduces the well-posedness of the mild solution to \eqref{eq14.1.1}.

\begin{lemma}\cite[Proposition 3.3]{DeVecchi2026}\label{l14.2.3}
	Assume (A1) and (A2). Then \eqref{eq14.1.1} admits a unique mild solution adapted to $W(t)$ taking values in $\mathcal{C}([0,T];H)$ of the form 
	\begin{equation*}
		X(t) = S(t)\xi + \int_0^tS(t-s)F(X(s))ds + \int_0^tS(t-s)A^{\gamma}dW(s).
	\end{equation*}
\end{lemma}

Define a stochastic convolution by 
\begin{equation*}
	\mathcal{O}_t = \int_0^tS(t-s)A^{\gamma}dW(s) = \int_0^te^{-A(t-s)}A^{\gamma}dW(s).
\end{equation*}
By means of the Kolmogorov test \cite[Section 5.5.1]{DaPrato2014}, $\mathcal{O}_t\in\mathcal{C}([0,T];\mathcal{C}(\bar{\mathcal{D}}))$ $a.s.$ 
In what follows, we employ a Da Prato-Debussche trick \cite{DaPrato2007} to investigate a priori estimates for the mild solution to \eqref{eq14.1.1}.
For any $\varrho\ge 0$, define a modified Ornstein-Uhlenbeck process by
\begin{equation*}
	\mathcal{O}_t^{\varrho} = \int_0^te^{-(A+\varrho)(t-s)}A^{\gamma}dW(s),
\end{equation*}
which is a mild solution to an SPDE
\begin{equation*}
	d\mathcal{O}_t^{\varrho} + (A+\varrho)\mathcal{O}_t^{\varrho}dt = A^{\gamma}dW(t),\quad \mathcal{O}_0^{\varrho} = 0.
\end{equation*}
Also by the Kolmogorov test, $\mathcal{O}^{\varrho}\in\mathcal{C}([0,T];\mathcal{C}(\bar{\mathcal{D}})) \subset \mathcal{C}([0,T];L^{\infty}(\mathcal{D}))$ $a.s.$ Let $\beta_j(t)$ be a family of standard $i.i.d.$ Wiener process, then we have
\begin{equation*}
	\mathcal{O}_t^{\varrho} = \sum_{j=1}^\infty \int_0^te^{-(\lambda_j+\varrho)(t-s)}\lambda_j^{\gamma}\varphi_jd\beta_j(s).
\end{equation*}
The following lemma collects some properties of $\mathcal{O}_t^{\varrho}$.

\begin{lemma}\label{l14.2.2}\cite[Proposition 3.7]{DeVecchi2026}
	Let $p>1$, $\varepsilon>0$, $\kappa\ge 0$ satisfying $\varepsilon+\kappa<\frac{1}{4}-\gamma$, and $\varrho\ge1$. Then $\mathcal{O}_t^{\varrho}$ take values in $\mathcal{C}([0,T];W^{2\kappa,p}\cap\mathcal{C}(\bar{\mathcal{D}}))$ $a.s.$ and is continuous with respect to $\varrho$. Moreover, there exists some random variables $K_{\varepsilon,p}^{(\kappa)}:\Omega\to\mathbb{R}_+$ which have finite moments of all orders, such that
	\begin{equation*}
		\|A^{\kappa}\mathcal{O}_t^{\varrho}\|_{L^p(\mathcal{D})} \le \varrho^{-\frac{1}{4}+\varepsilon+\gamma+\kappa}C_TK_{\varepsilon,p}^{(\kappa)},\quad a.s.
	\end{equation*}
\end{lemma}

Take $\varepsilon,\kappa> 0$ such that $\varepsilon+\kappa<\frac{1}{4}-\gamma$ and choose $p_0=\frac{1}{2\kappa}+1>2$. Then $2\kappa p_0>1$. According to the Sobolev embedding inequality $W^{2\kappa,p_0}\hookrightarrow V$, there exists a constant $C=C(\kappa)>0$ such that
\begin{equation}\label{eq14.2.9.1}
	\sup_{t\in[0,T]}\|\mathcal{O}_t^{\varrho}\|_V \le C\sup_{t\in[0,T]}\|A^{\kappa}\mathcal{O}_t^{\varrho}\|_{L^{p_0}(\mathcal{D})} \le C\varrho^{-\frac{1}{4}+\varepsilon+\gamma+\kappa}C_TK_{\varepsilon,p_0}^{(\kappa)},\quad a.s.
\end{equation}
Notice that $-\frac{1}{4}+\varepsilon+\gamma+\kappa<0$ and $K_{\varepsilon,p_0}^{(\kappa)}$ has finite moments of all orders, we obtain for $q\ge1$,
\begin{equation}\label{eq14.2.9}
	\mathbb{E}\sup_{t\in[0,T]}\|\mathcal{O}_t^{\varrho}\|_V^q \le C^q\varrho^{q(-\frac{1}{4}+\varepsilon+\gamma+\kappa)}C_T^q\mathbb{E}|K_{\varepsilon,p_0}^{(\kappa)}|^q \le C(T,q,\varepsilon,\kappa).
\end{equation}

The next theorem investigates a priori estimates for the mild solution $X(t)$ to \eqref{eq14.1.1}.

\begin{theorem}\label{t14.2.6}
	Assume (A1) and (A2). Let $X(t)$ be the mild solution to \eqref{eq14.1.1}, then for any $p,q\ge2$, there exists a constant $C>0$ such that for all $t\in[0,T]$,
		\begin{align*}
			\|X(t)\|_{L^q(\Omega;L^p(\mathcal{D}))} \le&\ C(1+\|\xi\|_{L^p(\mathcal{D})}),\\
			\mathbb{E}\sup_{t\in[0,T]}\|X(t)\|_V^q \le&\ C(1+\|\xi\|^2_{V})^q.
	\end{align*}
\end{theorem}

\begin{proof}
	Define an auxiliary process $Y^{\varrho}(t) = X(t)-\mathcal{O}_t^{\varrho}$, then $Y^{\varrho}(t)$ satisfies $Y^{\varrho}(0)=\xi$ and
	\begin{equation}\label{eq14.2.6}
		\frac{d}{dt}Y^{\varrho}(t,x) + AY^{\varrho}(t,x) = F(Y^{\varrho}(t,x)+\mathcal{O}_t^{\varrho}(x)) + \varrho\mathcal{O}_t^{\varrho}(x).
	\end{equation}
	Note that
	\begin{equation}\label{eq14.2.7}
		\nabla(|Y^{\varrho}(t,x)|^{p}Y^{\varrho}(t,x)) = (p+1)|Y^{\varrho}(t,x)|^{p}\nabla Y^{\varrho}(t,x).
	\end{equation}
	Testing \eqref{eq14.2.6} by $|Y^{\varrho}(t,x)|^{p-2}Y^{\varrho}(t,x)$, by the integration by parts formula and \eqref{eq14.2.7}, we obtain
	\begin{align*}
		&\ \frac{1}{p}\frac{d}{dt}\int_{\mathcal{D}}|Y^{\varrho}(t,x)|^{p}dx + \int_{\mathcal{D}}|Y^{\varrho}(t,x)|^{p-2}Y^{\varrho}(t,x)AY^{\varrho}(t,x)dx\\
		=&\ \int_{\mathcal{D}}|Y^{\varrho}(t,x)|^{p-2}Y^{\varrho}(t,x)\nabla[(Y^{\varrho}(t,x)+\mathcal{O}_t^{\varrho}(x))^2]dx + \varrho\int_{\mathcal{D}}|Y^{\varrho}(t,x)|^{p-2}Y^{\varrho}(t,x)\mathcal{O}_t^{\varrho}(x)dx\\
		=&\ -\int_{\mathcal{D}}\nabla[|Y^{\varrho}(t,x)|^{p-2}Y^{\varrho}(t,x)](Y^{\varrho}(t,x)+\mathcal{O}_t^{\varrho}(x))^2dx + \varrho\int_{\mathcal{D}}|Y^{\varrho}(t,x)|^{p-2}Y^{\varrho}(t,x)\mathcal{O}_t^{\varrho}(x)dx\\
		=&\ -(p-1)\int_{\mathcal{D}}|Y^{\varrho}(t,x)|^{p-2}\nabla Y^{\varrho}(t,x)(Y^{\varrho}(t,x)+\mathcal{O}_t^{\varrho}(x))^2dx + \varrho\int_{\mathcal{D}}|Y^{\varrho}(t,x)|^{p-2}Y^{\varrho}(t,x)\mathcal{O}_t^{\varrho}(x)dx\\
		=&\ -(p-1)\int_{\mathcal{D}}|Y^{\varrho}(t,x)|^{p}\nabla Y^{\varrho}(t,x))dx -2(p-1)\int_{\mathcal{D}}|Y^{\varrho}(t,x)|^{p-2}\nabla Y^{\varrho}(t,x) Y^{\varrho}(t,x)\mathcal{O}_t^{\varrho}(x)dx\\ 
		&\ -(p-1)\int_{\mathcal{D}}|Y^{\varrho}(t,x)|^{p-2}\nabla Y^{\varrho}(t,x)(\mathcal{O}_t^{\varrho}(x))^2dx + \varrho\int_{\mathcal{D}}|Y^{\varrho}(t,x)|^{p-2}Y^{\varrho}(t,x)\mathcal{O}_t^{\varrho}(x)dx\\
		 =:&\ I_1+I_2+I_3+I_4.
	\end{align*}
	
	We first calculate the second term in the left hand side of above equation
	\begin{align*}
		&\ \int_{\mathcal{D}}AY^{\varrho}(t,x)|Y^{\varrho}(t,x)|^{p-2}Y^{\varrho}(t,x)dx = \int_{\mathcal{D}}\nabla_x Y^{\varrho}(t,x)\nabla_x[|Y^{\varrho}(t,x)|^{p-2}Y^{\varrho}(t,x)]dx\\ 
		=&\ (p-1)\int_{\mathcal{D}}|Y^{\varrho}(t,x)|^{p-2}(\nabla_xY^{\varrho}(t,x))^2dx = (p-1)\||Y^{\varrho}(t)|^{{\frac{p}{2}}-1}\nabla Y^{\varrho}(t)\|^{2},
	\end{align*}
	
	Then we study the four terms in the right hand side.  By \eqref{eq14.2.7}, we have
	\begin{equation*}
		I_1 = -\frac{p-1}{p+1}\int_{0}^1\nabla(|Y^{\varrho}(t,x)|^{p}Y^{\varrho}(t,x))dx = 0.
	\end{equation*}
	According to \cite{DaPrato2007}, we have $\nabla|Y^{\varrho}(t)|^{\frac{p}{2}} = \frac{p}{2}|Y^{\varrho}(t)|^{\frac{p}{2}-2}Y^{\varrho}(t)\nabla Y^{\varrho}(t) = \frac{p}{2}\mathrm{sgn}(Y^{\varrho}(t))|Y^{\varrho}(t)|^{\frac{p}{2}-1}\nabla Y^{\varrho}(t)$. By Sobolev embedding $\dot{H}^{\frac{1}{4}}\hookrightarrow L^4(\mathcal{D})$ and the Gagliardo-Nirenberg inequality, we obtain
	\begin{align*}
		&\ \||Y^{\varrho}(t)|^{\frac{p}{2}}\|_{L^4(\mathcal{D})} \le C\||Y^{\varrho}(t)|^{\frac{p}{2}}\|_{\frac{1}{4}} \le C\||Y^{\varrho}(t)|^{\frac{p}{2}}\|^{\frac{3}{4}}\cdot\||Y^{\varrho}(t)|^{\frac{p}{2}}\|^{\frac{1}{4}}_1\\
		\le&\ C\||Y^{\varrho}(t)|^{\frac{p}{2}}\|^{\frac{3}{4}}\cdot C\|\nabla|Y^{\varrho}(t)|^{\frac{p}{2}}\|^{\frac{1}{4}} = C\|Y^{\varrho}(t)\|_{L^p(\mathcal{D})}^{\frac{3}{8}p}\cdot \||Y^{\varrho}(t)|^{\frac{p}{2}-1}\nabla Y^{\varrho}(t)\|^{\frac{1}{4}},
	\end{align*}
	thus by the H\"older inequality and Young's inequality $ab\le \frac{5}{8}a^{\frac{8}{5}} + \frac{3}{8}b^{\frac{8}{3}}$, we get
	\begin{align*}
		|I_2| \le&\ C\int_{\mathcal{D}}|Y^{\varrho}(t,x)|^{\frac{p}{2}}|Y^{\varrho}(t,x)|^{\frac{p}{2}-1}|\nabla Y^{\varrho}(t,x)|\cdot|\mathcal{O}_t^{\varrho}(x)|dx\\
		\le&\  C\||Y^{\varrho}(t)|^{\frac{p}{2}}\|_{L^4(\mathcal{D})}\cdot \||Y^{\varrho}(t)|^{{\frac{p}{2}}-1}\nabla Y^{\varrho}(t)\|\cdot \|\mathcal{O}_t^{\varrho}\|_{L^4(\mathcal{D})}\\ 
		\le&\ C\|Y^{\varrho}(t)\|_{L^p(\mathcal{D})}^{\frac{3}{8}p}\||Y^{\varrho}(t)|^{{\frac{p}{2}}-1}\nabla Y^{\varrho}(t)\|^{\frac{1}{4}}\cdot \||Y^{\varrho}(t)|^{{\frac{p}{2}}-1}\nabla Y^{\varrho}(t)\|\cdot\|\mathcal{O}_t^{\varrho}\|_{L^4(\mathcal{D})}\\
		=&\ \Bigl(\frac{2(p-1)}{5}\Bigr)^{\frac{5}{8}}\||Y^{\varrho}(t)|^{\frac{p}{2}-1}\nabla Y^{\varrho}(t)\|^{\frac{5}{4}}\cdot C\Bigl(\frac{5}{2(p-1)}\Bigr)^{\frac{5}{8}}\|Y^{\varrho}(t)\|_{L^p(\mathcal{D})}^{\frac{3p}{8}}\|\mathcal{O}_t^{\varrho}\|_{L^4(\mathcal{D})}\\
		\le&\ \frac{p-1}{4}\||Y^{\varrho}(t)|^{\frac{p}{2}-1}\nabla Y^{\varrho}(t)\|^2+\frac{3}{8}C^{\frac{8}{3}}\Bigl(\frac{5}{2(p-1)}\Bigr)^{\frac{5}{3}}\|Y^{\varrho}(t)\|_{L^p(\mathcal{D})}^p\|\mathcal{O}_t^{\varrho}\|_{L^4(\mathcal{D})}^{\frac{8}{3}}.
	\end{align*}
	By the Gagliardo-Nirenberg inequality, we have
	\begin{equation*}
		\|Y^{\varrho}(t)^{{\frac{p}{2}}}\|_V \le \|Y^{\varrho}(t)^{{\frac{p}{2}}}\|^{\frac{1}{2}}\|Y^{\varrho}(t)^{{\frac{p}{2}}}\|^{\frac{1}{2}}_1 \le C\|Y^{\varrho}(t)\|_{L^p(\mathcal{D})}^{\frac{p}{4}}\||Y^{\varrho}(t)|^{{\frac{p}{2}}-1}\nabla Y^{\varrho}(t)\|^{\frac{1}{2}},
	\end{equation*}
	which implies
	\begin{equation*}
		\|Y^{\varrho}(t)\|_V^{\frac{p}{2}-1}= \bigl(\|Y^{\varrho}(t)\|_V^{\frac{p}{2}}\bigr)^{\frac{p-2}{p}} = \bigl(\|Y^{\varrho}(t)^{\frac{p}{2}}\|_V\bigr)^{\frac{p-2}{p}} \le \|Y^{\varrho}(t)\|_{L^p(\mathcal{D})}^{\frac{p-2}{4}}\||Y^{\varrho}(t)|^{{\frac{p}{2}}-1}\nabla Y^{\varrho}(t)\|^{\frac{p-2}{2p}}.
	\end{equation*}
	Thus by Young's inequality $abc\le \frac{3p-2}{4p}a^{\frac{4p}{3p-2}} + \frac{p-2}{4p}b^{\frac{4p}{p-2}} + \frac{1}{p}c^p$, we obtain
	\begin{align*}
		|I_3| \le&\  C\int_{\mathcal{D}}|Y^{\varrho}(t,x)|^{\frac{p}{2}-1}|\nabla Y^{\varrho}(t,x)|\cdot|Y^{\varrho}(t,x)|^{\frac{p}{2}-1}\cdot(\mathcal{O}_t^{\varrho}(x))^2dx\\
		\le&\ C\||Y^{\varrho}(t)|^{\frac{p}{2}-1}\nabla Y^{\varrho}(t)\|\|Y^{\varrho}(t)\|_V^{\frac{p}{2}-1}\|\mathcal{O}_t^{\varrho}\|_{L^4(\mathcal{D})}^2\\
		\le&\ C\||Y^{\varrho}(t)|^{\frac{p}{2}-1}\nabla Y^{\varrho}(t)\|\cdot\|Y^{\varrho}(t)\|_{L^p(\mathcal{D})}^{\frac{p-2}{4}}\||Y^{\varrho}(t)|^{{\frac{p}{2}}-1}\nabla Y^{\varrho}(t)\|^{\frac{p-2}{2p}}\cdot\|\mathcal{O}_t^{\varrho}\|_{L^4(\mathcal{D})}^2\\
		\le&\ C\||Y^{\varrho}(t)|^{{\frac{p}{2}}-1}\nabla Y^{\varrho}(t)\|^{\frac{3p-2}{2p}}\|Y^{\varrho}(t)\|_{L^p(\mathcal{D})}^{\frac{p-2}{4}}\|\mathcal{O}_t^{\varrho}\|_{L^4(\mathcal{D})}^2\\
		=&\ \Bigl(\frac{p(p-1)}{3p-2}\Bigr)^{\frac{3p-2}{4p}}\||Y^{\varrho}(t)|^{{\frac{p}{2}}-1}\nabla Y^{\varrho}(t)\|^{\frac{3p-2}{2p}}\cdot \|Y^{\varrho}(t)\|_{L^p(\mathcal{D})}^{\frac{p-2}{4}}\cdot C\Bigl(\frac{3p-2}{p(p-1)}\Bigr)^{\frac{3p-2}{4p}}\|\mathcal{O}_t^{\varrho}\|_{L^4(\mathcal{D})}^2\\
		\le&\ \frac{p-1}{4}\||Y^{\varrho}(t)|^{\frac{p}{2}-1}\nabla Y^{\varrho}(t)\|^2+C(p)\|Y^{\varrho}(t)\|_{L^p(\mathcal{D})}^p+C(p)\|\mathcal{O}_t^{\varrho}\|_{L^4(\mathcal{D})}^{2p}.
	\end{align*}
	For $p=2$, Young's inequality yields
	\begin{equation*}
		|I_3|\le\|\nabla Y^{\varrho}(t)\|\|\mathcal{O}_t^{\varrho}\|_{L^4(\mathcal{D})}^2\le\frac{1}{4}\|\nabla Y^{\varrho}(t)\|^2+\|\mathcal{O}_t^{\varrho}\|_{L^4(\mathcal{D})}^4.
	\end{equation*}
	By the H\"older inequality and Young's inequality, we have
	\begin{equation*}
		|I_4| \le \varrho\|Y^{\varrho}(t)\|_{L^p(\mathcal{D})}^{p-1}\|\mathcal{O}_t^{\varrho}\|_{L^p(\mathcal{D})} \le \frac{p-1}{p}\|Y^{\varrho}(t)\|_{L^p(\mathcal{D})}^{p} + \frac{1}{p}\varrho^p\|\mathcal{O}_t^{\varrho}\|_{L^p(\mathcal{D})}^p.
	\end{equation*}
	Above all, we obtain
	\begin{align*}
		&\ \frac{d}{dt}\|Y^{\varrho}(t)\|_{L^p(\mathcal{D})}^{p} + p(p-1)\||Y^{\varrho}(t)|^{{\frac{p}{2}}-1}\nabla Y^{\varrho}(t)\|^{2} \le \frac{p(p-1)}{2}\||Y^{\varrho}(t)|^{{\frac{p}{2}}-1}\nabla Y^{\varrho}(t)\|^{2}\\ 
		&\ + C\Bigl(\|\mathcal{O}_t^{\varrho}\|_{L^4(\mathcal{D})}^{\frac{8}{3}} + \|\mathcal{O}_t^{\varrho}\|_{L^p(\mathcal{D})} + 1\Bigr)\|Y^{\varrho}(t)\|_{L^p(\mathcal{D})}^p + C\bigl(\varrho^p\|\mathcal{O}_t^{\varrho}\|^p_{L^p(\mathcal{D})} + \|\mathcal{O}_t^{\varrho}\|_{L^4(\mathcal{D})}^{2p}\bigr),
	\end{align*}
	which leads to
	\begin{equation*}
		\frac{d}{dt}\|Y^{\varrho}(t)\|_{L^p(\mathcal{D})}^{p} \le C\Bigl(\|\mathcal{O}_t^{\varrho}\|_{L^4(\mathcal{D})}^{\frac{8}{3}} + \|\mathcal{O}_t^{\varrho}\|_{L^p(\mathcal{D})} + 1\Bigr)\|Y^{\varrho}(t)\|_{L^p(\mathcal{D})}^p + C\bigl(\varrho^p\|\mathcal{O}_t^{\varrho}\|^p_{L^p(\mathcal{D})} + \|\mathcal{O}_t^{\varrho}\|_{L^4(\mathcal{D})}^{2p}\bigr).
	\end{equation*}
	By Lemma \ref{l14.2.2}, we choose $0<\varepsilon<\frac{1}{4}-\gamma$, $\kappa=0$, $K=\max\{K^{(0)}_{\varepsilon,4},K^{(0)}_{\varepsilon,p}\}$, and $\varrho = \bigl(C_TK\bigr)^{\frac{1}{1/4-\varepsilon-\gamma}}+1$, then $\varrho^{-\frac{1}{4}+\varepsilon+\gamma}C_TK\le1$ and
	\begin{equation*}
		\sup_{t\in[0,T]}\|\mathcal{O}_t^{\varrho}\|_{L^4(\mathcal{D})}^{\frac{8}{3}} + \sup_{t\in[0,T]}\|\mathcal{O}_t^{\varrho}\|_{L^p(\mathcal{D})} \le \varrho^{-\frac{1}{4}+\varepsilon+\gamma}C_TK^{(0)}_{\varepsilon,4} + \varrho^{-\frac{1}{4}+\varepsilon+\gamma}C_TK^{(0)}_{\varepsilon,p} \le 1+1=2,\ a.s.
	\end{equation*}
	Then we get
	\begin{equation*}
		\frac{d}{dt}\|Y^{\varrho}(t)\|_{L^p(\mathcal{D})}^{p} \le C\|Y^{\varrho}(t)\|_{L^p(\mathcal{D})}^p + C\bigl(\varrho^p + 1\bigr),\quad a.s.
	\end{equation*}
	From a Gronwall's inequality, it follows that
	\begin{equation*}
		\|Y^{\varrho}(t)\|_{L^p(\mathcal{D})}^{p} \le \bigl(\|Y^{\varrho}(0)\|^p_{L^p(\mathcal{D})}+Ct\bigl(\varrho^p + 1\bigr)\bigr)e^{Ct} \le C\bigl(\|\xi\|_{L^p}^p+1+\varrho(\omega)^p\bigr).,\quad a.s.
	\end{equation*}
	Thus by Lemma \ref{l14.2.3}, we obtain for any $q\ge 2$,
	\begin{equation}\label{eq14.2.8}
		\|X(t)\|_{L^q(\Omega;L^p(\mathcal{D}))} \le \|Y^{\varrho}(t)\|_{L^q(\Omega;L^p(\mathcal{D}))} + \|\mathcal{O}_t^\varrho\|_{L^q(\Omega;L^p(\mathcal{D}))} \le C(1+\|\xi\|_{L^p(\mathcal{D})}).
	\end{equation}
	
	Now we estimate $X(t)$ in $V$-norm. According to \cite[Lemma 3.3]{CuiJianbo2019}, we have 
	\begin{equation*}
		\sup_{t\in[0,T]}\|S(t)\xi\|_{V} \le C\|\xi\|_V.
	\end{equation*} 
	Note that
	\begin{equation*}
		\mathcal O_t = \mathcal O_t^\varrho + \varrho\int_0^tS(t-s)\mathcal O_s^\varrho\,ds,
	\end{equation*}
	In particular, taking $\varrho=1$, we obtain
	\begin{align*}
		\sup_{t\in[0,T]}\|\mathcal O_t\|_V
		&\le
		\sup_{t\in[0,T]}\|\mathcal O_t^1\|_V
		+
		\sup_{t\in[0,T]}
		\int_0^t\|S(t-s)\mathcal O_s^1\|_V\,ds\\
		&\le
		(1+T)\sup_{t\in[0,T]}\|\mathcal O_t^1\|_V.
	\end{align*}
	Consequently, by \eqref{eq14.2.9},
	\begin{equation*}
		\mathbb E\sup_{t\in[0,T]}\|\mathcal O_t\|_V^q \le C_{T,q} \mathbb E\sup_{t\in[0,T]}\|\mathcal O_t^1\|_V^q \le C.
	\end{equation*}
	For small $0<\delta<\frac{1}{2}$, take $m> \max\{\frac{4}{1-2\delta},q\}$, then $(\frac{3}{4}+\frac{\delta}{2})\frac{m}{m-1}<1$ and $\frac{q}{m}<1$. By \eqref{eq12.2.3}, Sobolev embedding property $\dot{H}^{\frac{1}{2}+\delta}\hookrightarrow V$, and the H\"older inequality, we have
	\begin{align*}
		&\ \Bigl\|\int_0^tS(t-s)F(X(s))ds\Bigr\|_V \le \int_0^t\|A^{\frac{3}{4}+\frac{\delta}{2}}S(t-s)A^{-\frac{1}{2}}\partial_xX^2(s)\|ds\\
		\le&\ C\int_0^t(t-s)^{-(\frac{3}{4}+\frac{\delta}{2})}\|X(s)\|_{L^4(\mathcal{D})}^2ds \le C\Bigl(\int_0^t(t-s)^{-(\frac{3}{4}+\frac{\delta}{2})\frac{m}{m-1}}ds\Bigr)^{\frac{m-1}{m}}\Bigl(\int_0^t\|X(s)\|_{L^4(\mathcal{D})}^{2m}\Bigr)^{\frac{1}{m}}\\
		\le&\ Ct^{\frac{1}{4}-\frac{\delta}{2}-\frac{1}{m}}\Bigl(\int_0^t\|X(s)\|_{L^4(\mathcal{D})}^{2m}\Bigr)^{\frac{1}{m}}.
	\end{align*}
	Since $\frac{q}{m}<1$ and Jensen's inequality, we have $\mathbb{E}[\|X\|^{\frac{q}{m}}]\le (\mathbb{E}\|X\|)^{\frac{q}{m}}$ for any random variable $X$. Note that $\frac{1}{4}-\frac{\delta}{2}-\frac{1}{m}>0$, we have
	\begin{align*}
		\mathbb{E}\sup_{t\in[0,T]}\|X(t)\|_V^q\le&\ C\sup_{t\in[0,T]}\|S(t)\xi\|_V^q+C\mathbb{E}\sup_{t\in[0,T]}\Bigl\|\int_0^tS(t-s)F(X(s))ds\Bigr\|_V^q+C\mathbb{E}\sup_{t\in[0,T]}\|\mathcal{O}_t\|_V^q\\
		\le&\ C\|\xi\|_V^q+CT^{q(\frac{1}{4}-\frac{\delta}{2}-\frac{1}{m})}\mathbb{E}\Bigl(\int_0^T\|X(s)\|_{L^4(\mathcal{D})}^{2m}ds\Bigr)^{\frac{q}{m}}+C\\
		\le&\ C(1+\|\xi\|_V^q)+C\Bigl(\int_0^T\mathbb{E}\|X(s)\|_{L^4(\mathcal{D})}^{2m}ds\Bigr)^{\frac{q}{m}}\\
		\le&\ C(1+\|\xi\|_V^q)+C(1+\|\xi\|_{L^4(\mathcal{D})})^{2q}\le C(1+\|\xi\|_V^2)^q.
	\end{align*}
	Thus the proof is completed.
\end{proof}

\subsection{Spectral Galerkin approximation}

For any integer $N\ge 1$, define a finite dimensional Galerkin subspace $H_N = \mathrm{Span}\{\varphi_1,\cdots,\varphi_N\}$. Let $P_N:H \to H_N$ be an orthogonal projection operator,  then $P_NA=AP_N$, \cite{CaiMeng2023}. According to \cite[Lemma 4.1]{QiRuisheng2019}, there hold
\begin{equation}\label{eq14.5.1}
	 \|(P_N-I)A^{-\alpha}\|_{\mathcal{L}(H)}\le CN^{-2\alpha},\quad \forall\ \alpha\in[0,\infty),
\end{equation}
where $C=C(\alpha)>0$ is a constant independent of $N$.
The spectral Galerkin approximation to \eqref{eq14.1.1} is defined as: seeking $X^{N}(t)\in H_N$ such that
\begin{equation}\label{eq14.5.3}
	\begin{aligned}
		dX^{N}(t) + AX^{N}(t)dt =&\ P_NF(X^{N}(t))dt + P_NA^{\gamma}dW(t),\ t\in (0,T],\\
		X^{N}(0) =&\ P_N \xi.
	\end{aligned}
\end{equation}
Its mild solution is given by
\begin{equation}\label{eq14.5.4}
	X^{N}(t) = S(t)P_N\xi + \int_0^tS(t-s)P_NF(X^{N}(s))ds + \int_0^tS(t-s)P_NA^{\gamma}dW(s),\quad \forall t\in[0,T].
\end{equation}
We can apply a similar procedure as proving Theorem \ref{t14.2.6} to derive a priori estimates for $X^N(t)$, which is summarized in next lemma and the proof is omitted.

\begin{lemma}\label{l14.5.1}
	Assume (A1) and (A2). Then for any $p,q\ge2$, there exists a constant $C=C(T)>0$ independent of $N$ such that for all $t\in[0,T]$,
	\begin{align*}
		\|X^N(t)\|_{L^q(\Omega;L^p(\mathcal{D}))} \le&\ C(1+\|\xi\|_{L^p(\mathcal{D})}),\\
		\mathbb{E}\sup_{t\in[0,T]}\|X^N(t)\|_V^q \le&\ C(1+\|P_N\xi\|^2_{V})^q.
	\end{align*}
\end{lemma}

We would like to emphasize that (A2) ensures the projected initial data are uniformly bounded in the maximum norm.
In fact, due to the Sobolev embedding inequality $\dot{H}^1\hookrightarrow V$ and $\|P_Nu\|_1\le \|u\|_1$ for all $u\in \dot{H}^1$, there holds
\begin{equation*}
	\|P_N\xi\|_{V} \le C\|A^{\frac{1}{2}}P_N\xi\| = C\|P_NA^{\frac{1}{2}}\xi\| \le C\|A^{\frac{1}{2}}\xi\| = C\|\xi\|_1 <\infty.
\end{equation*}
This guarantees the uniformly boundedness of $\mathbb{E}\sup_{t\in[0,T]}\|X^N(t)\|_V^q$.

The following lemma describes the convergence of the Galerkin approximation $X^N$ to $X$ as $N\to\infty$.

\begin{lemma}\label{l14.3.1}\cite[Theorem 3.4]{DeVecchi2026}
	Assume (A1) and (A2). Then the sequence $X^{N}$ converges to $X$ in $\mathcal{C}([0,T];H)$ as $N\to \infty$ $a.s.$ 
\end{lemma}

Denote by $\mathcal{P}^{[0,T],H_N}$ the space of $H_N$-valued predictable process with bounded second moment.
For any $u\in\mathcal{P}^{[0,T],H_N}$, let $X^{N,u}(t)\in H_N$ be the solution of the following auxiliary SDE
\begin{equation}\label{eq14.3.2}
	\begin{aligned}
		dX^{N,u}(t)+A X^{N,u}(t)dt =&\ P_NF(X^{N,u}(t))dt + P_NA^{\gamma}u(t)dt +P_NA^{\gamma}dW(t),\\
		X^{N,u}(0) =&\ x_N.
	\end{aligned}
\end{equation}
Define two auxiliary processes in $H_N$ by $Y^{N,\varrho}(t) = X^{N}(t) - P_N\mathcal{O}_t^{\varrho}$ and $Y^{N,u,\varrho}(t) = X^{N,u}(t) - P_N\mathcal{O}_t^{\varrho}$. Then as $u=0$, we have $Y^{N,0,\varrho}(t)=Y^{N,\varrho}(t)$. The following lemma presents a priori estimates of the solution to \eqref{eq14.3.2}.

\begin{lemma}\label{l14.3.2}\cite[Theorem 3.10]{DeVecchi2026}
	Assume (A1) and (A2). Then for any $\delta_1>0$, $0<\varepsilon_1<\frac{1}{4}-\gamma$, and $\theta:=2(\frac{1}{4}-\varepsilon_1-\gamma)^{-1}>0$, there exists a constant $C_{T,\varepsilon_1}>0$ such that for any $t\in[0,T]$, there holds $a.s.$
	\begin{equation*}
		\|X^{N,u}(t)\|^2 \le e^{-\delta_1t}\|x_N\|^2 + C_{T,\varepsilon_1}\bigl(1+K_{\varepsilon_1,4}^{\theta}+\|u\|^2_{L^2_tL^2_{x}}\bigr),
	\end{equation*}
	where $K_{\varepsilon_1,4}$ is given in Lemma \ref{l14.2.2}.
	Moreover,  for any small parameter $0<\varepsilon'\le 1$, square integrable random variable $G$, and a random variable $\varrho_0$ of the form
	\begin{equation*}
		\varrho_0(\omega) = \frac{1}{\varepsilon'}\Bigl(C_TK^{\frac{\theta}{2}}_{\varepsilon_1,4}(\omega)+1+G(\omega)\Bigr),
	\end{equation*}
	there exists constants $C_{\gamma,\varepsilon_1},\tilde{C}_{T,\varepsilon'}>0$ such that for any $t\in[0,T]$, there holds $a.s.$
	\begin{align*}
		\|Y^{N,u,\varrho_0(\omega)}(t)\|^2 \le&\ e^{-\delta_1t}\|x_N\|^2 + \tilde{C}_{T,\varepsilon'}\bigl(1+ G^2(\omega)+K_{\varepsilon_1,4}^{\theta}+\|u\|^2_{L^2_tL^2_{x}}\bigr),\\
		\int_0^T\|\nabla Y^{N,u,\varrho_0(\omega)}(s)\|^2ds \le&\ C_{\gamma,\varepsilon_1}\|x_N\|^2 + \tilde{C}_{T,\varepsilon'}\bigl(1+ G^2(\omega) + K_{\varepsilon_1,4}^{\theta}+\|u\|^2_{L^2_tL^2_{x}}\bigr).
	\end{align*}
\end{lemma}

For any given square integrable random variable $G$, let $\varrho_0$ be defined in Lemma \ref{l14.3.2}. The next lemma states an exponential integrability of the exact solution $X(t)$ and the deterministic part $Y^{N,\varrho_0}(t)$ of the spectral Galerkin approximation $X^N(t)$.

\begin{lemma}\label{t14.3.1}\cite[Theorem 3.12]{DeVecchi2026}
	Assume (A1) and (A2). Then there exists constants $\lambda_{0}>0$ and $C_{\gamma,T}>0$ such that for all $0\le\lambda\le\lambda_{0}$, there holds
	\begin{equation*}
		\mathbb{E}\Bigl[\exp\Bigl(\lambda\sup_{t\in[0,T]}\|X(t)\|^2\Bigr)\Bigr] + \mathbb{E}\Bigl[\exp\Bigl(\lambda\int_0^T\|\nabla Y^{N,\varrho_0}(t)\|^2dt\Bigr)\Bigr] \le C_{\gamma,T}e^{\lambda\|\xi\|^2}.
	\end{equation*}
\end{lemma}

\section{Kolmogorov equation and regularity estimates}
\label{Sect3}

For any given $\nu\in H_N$, let $X^N(t,\nu)$ be the solution of \eqref{eq14.5.3} with initial value $X^{N}(0)=\nu$.
For any $\phi \in \mathcal{C}_b^2(H)$, define $U^N(t,\nu) = \mathbb{E}[\phi(X^N(t,\nu))]$, which satisfies a Kolmogorov equation with initial value $U^N(0,\nu)=\phi(\nu)$,
\begin{equation}\label{eq14.4.1}
	\frac{\partial U^N(t,\nu)}{\partial t}= \bigl(DU^N(t,\nu),-A\nu+P_NF(\nu)\bigr) +\frac{1}{2}\mathrm{Tr}\bigl\{D^2U^N(t,\nu)P_NA^{2\gamma}P_N^*\bigr\}.
\end{equation}
The main objective of this section is to investigate a priori estimates for the first and second order Fr\'echet derivatives $DU^N(t,\nu)$ and $D^2U^N(t,\nu)$ of $U^N(t,\nu)$ with respect to $\nu$, including estimates involving fractional powers of Laplacian operator, which will provide a weak convergence order to the spectral Galerkin approximation.

Let $\eta^h(t) = DX^N(t,\nu)h$ be the directional Fr\'echet derivative of $X^N(t,\nu)$ with respect to $\nu$ along $h\in H_N$.
Taking directional Fr\'echet derivative on both sides of \eqref{eq14.5.3} along $h\in H_N$, we obtain $\eta^h(0)=h$ and an equation that $\eta^h(t)$ satisfies
\begin{equation}\label{eq14.4.2}
	d\eta^h(t) + A\eta^h(t)dt = P_NDF(X^N(t))\eta^h(t)dt.
\end{equation}
Let $\Pi(t,s)$ with $\Pi(s,s)=I$ be the solution operator of \eqref{eq14.4.2}. Then, $\eta^h(t,s) = \Pi(t,s)h$ is a solution of \eqref{eq14.4.2} with $\eta^h(s,s)=h$, especially, $\eta^h(t)=\eta^h(t,0)$.

We first investigate pathwise a priori estimates for the solution operator $\Pi(t,s)$ to \eqref{eq14.4.2}.

\begin{lemma}\label{l14.4.2}
	Assume (A1) and (A2). Then for any $\alpha\in[0,1)$, $\varepsilon>0$, and $0<\rho<\frac{1}{16}$, there exists a constant $C > 0$ independent of $N$, such that for all $h \in H_N$ and $0\le s<t\le T$,
	\begin{align*}
		\|\Pi(t,s)A^{\alpha}h\| \le&\ C\sup_{u\in[0,T]}\|X^N(u)\|^2_V\int_s^t(r-s)^{-\alpha}\|h\|\exp\Bigl(\int_r^t \varepsilon\|\nabla Y^{N,\varrho}(u)\|^{2} + \frac{3}{\rho}\|P_N\mathcal{O}_u^{\varrho}\|^2_Vdu\Bigr)dr\\ 
		&\ +C(t-s)^{-\alpha}\sup_{u\in[0,T]}\|X^N(u)\|_V\|h\| + C(t-s)^{-\alpha}\|h\|,\quad a.s.
	\end{align*}
\end{lemma}

\begin{proof}
	For any $\alpha\in[0,1)$, define $\eta^{A^{\alpha}h}(t,s) = \Pi(t,s)A^{\alpha}h$ and $\tilde{\eta}^{A^{\alpha}h}(t,s) = \eta^{A^{\alpha}h}(t,s) - S(t-s)A^{\alpha}h$ with $\tilde{\eta}^{A^{\alpha}h}(s,s) = 0$. Substituting $\eta^{A^{\alpha}h}(t,s) = \tilde{\eta}^{A^{\alpha}h}(t,s) + S(t-s)A^{\alpha}h$ into \eqref{eq14.4.2} yields
	\begin{equation}\label{eq14.4.3}
		d\tilde{\eta}^{A^{\alpha}h}(t,s) + A\tilde{\eta}^{A^{\alpha}h}(t,s)dt = P_NDF(X^N(t))\tilde{\eta}^{A^{\alpha}h}(t,s)dt + P_NDF(X^N(t))S(t-s)A^{\alpha}hdt.
	\end{equation}
	
	We will derive pathwise estimates for $\tilde{\eta}^{A^{\alpha}h}(t,s)$ in the cases $\alpha\in[0,\frac{1}{2}]$ and $\alpha\in(\frac{1}{2},1)$, respectively.
	As $\alpha\in[0,\frac{1}{2}]$, testing \eqref{eq14.4.3} by $\tilde{\eta}^{A^{\alpha}h}(t,s)$ and integrating over $[s,t]$, we get
	\begin{align*}
		&\ \frac{1}{2}\|\tilde{\eta}^{A^{\alpha}h}(t,s)\|^2 + \int_s^t\bigl(A\tilde{\eta}^{A^{\alpha}h}(r,s),\tilde{\eta}^{A^{\alpha}h}(r,s)\bigr)dr\\
		=&\  \int_s^t\bigl(P_NDF(X^N(r))\tilde{\eta}^{A^{\alpha}h}(r,s),\tilde{\eta}^{A^{\alpha}h}(r,s)\bigr)dr + \int_s^t\bigl(P_NDF(X^N(r))S(r-s)A^{\alpha}h,\tilde{\eta}^{A^{\alpha}h}(r,s)\bigr)dr.
	\end{align*}
	First we note that
	\begin{equation*}
		\bigl(A\tilde{\eta}^{A^{\alpha}h}(r,s),\tilde{\eta}^{A^{\alpha}h}(r,s)\bigr) = \bigl(A^{\frac{1}{2}}\tilde{\eta}^{A^{\alpha}h}(r,s),A^{\frac{1}{2}}\tilde{\eta}^{A^{\alpha}h}(r,s)\bigr) = \|\tilde{\eta}^{A^{\alpha}h}(r,s)\|_1^2.
	\end{equation*}
	According to \eqref{eq14.2.5}, we get
	\begin{align*}
		&\ \bigl(P_NDF(X^N(r))\tilde{\eta}^{A^{\alpha}h}(r,s),\tilde{\eta}^{A^{\alpha}h}(r,s)\bigr) = \bigl(P_NDF(Y^{N,\varrho}(r)+P_N\mathcal{O}_r^{\varrho})\tilde{\eta}^{A^{\alpha}h}(r,s),\tilde{\eta}^{A^{\alpha}h}(r,s)\bigr)\\ 
		=&\  \Bigl(2\partial_x\big(Y^{N,\varrho}(r)+P_N\mathcal{O}_r^{\varrho}\big)\tilde{\eta}^{A^{\alpha}h}(r,s),\tilde{\eta}^{A^{\alpha}h}(r,s)\Bigr) + \Bigl(2\big(Y^{N,\varrho}(r)+P_N\mathcal{O}_r^{\varrho}\big)\partial_x\tilde{\eta}^{A^{\alpha}h}(r,s),\tilde{\eta}^{A^{\alpha}h}(r,s)\Bigr)\\
		=&\ 2\int_{\mathcal{D}}\partial_xY^{N,\varrho}(r)|\tilde{\eta}^{A^{\alpha}h}(r,s)|^2dx + 2\int_{\mathcal{D}}\partial_x(P_N\mathcal{O}_r^{\varrho})|\tilde{\eta}^{A^{\alpha}h}(r,s)|^2dx\\
		&\ + 2\int_{\mathcal{D}}Y^{N,\varrho}(r)\partial_x\tilde{\eta}^{A^{\alpha}h}(r,s) \tilde{\eta}^{A^{\alpha}h}(r,s)dx + 2\int_{\mathcal{D}}P_N\mathcal{O}_r^{\varrho}\partial_x\tilde{\eta}^{A^{\alpha}h}(r,s) \tilde{\eta}^{A^{\alpha}h}(r,s)dx\\
		=:&\ I+II+III+IV.
	\end{align*}
	The Gagliardo-Nirenberg inequality implies
	\begin{equation*}
		\|\tilde{\eta}^{A^{\alpha}h}(r,s)\|_{L^4(\mathcal{D})}^2 \le \sqrt{2}\|\tilde{\eta}^{A^{\alpha}h}(r,s)\|^{\frac{3}{2}}\|\tilde{\eta}^{A^{\alpha}h}(r,s)\|_1^{\frac{1}{2}}.
	\end{equation*}
	Then by the H\"older inequality, Young's inequality $ab\le \frac{3}{4}a^{\frac{4}{3}}+\frac{1}{4}b^4$, and $\varepsilon$-Young's inequality $3a^{\frac{4}{3}}\le \frac{16}{\varepsilon^2}+\frac{\varepsilon a^2}{2}$ with small $\varepsilon>0$, we have
	\begin{align*}
		I=&\ \bigl|2\int_{\mathcal{D}}\partial_xY^{N,\varrho}(r)|\tilde{\eta}^{A^{\alpha}h}(r,s)|^2dx\bigr| \le 2\|\nabla Y^{N,\varrho}(r)\|\cdot\|\tilde{\eta}^{A^{\alpha}h}(r,s)\|_{L^4(\mathcal{D})}^2\\
		\le&\ 2\sqrt{2}\|\nabla Y^{N,\varrho}(r)\|\cdot \|\tilde{\eta}^{A^{\alpha}h}(r,s)\|^{\frac{3}{2}}\|\tilde{\eta}^{A^{\alpha}h}(r,s)\|_1^{\frac{1}{2}}\\ 
		\le&\ 3\|\nabla Y^{N,\varrho}(r)\|^{\frac{4}{3}}\|\tilde{\eta}^{A^{\alpha}h}(r,s)\|^2 + \frac{1}{4}\|\tilde{\eta}^{A^{\alpha}h}(r,s)\|^2_1\\
		\le&\ \bigl(\frac{16}{\varepsilon^2}+\frac{\varepsilon}{2}\|\nabla Y^{N,\varrho}(r)\|^{2}\bigr)\|\tilde{\eta}^{A^{\alpha}h}(r,s)\|^2 + \frac{1}{4}\|\tilde{\eta}^{A^{\alpha}h}(r,s)\|^2_1.
	\end{align*}
	From integration by parts formula and $\rho$-Young's inequality $ab\le \frac{\rho^{-1}a^2}{2}+\frac{\rho b^2}{2}$ with $0<\rho<\frac{1}{16}$, it follows that
	\begin{align*}
		II=&\ \bigl|2\int_{\mathcal{D}}\partial_x(P_N\mathcal{O}_r^{\varrho})|\tilde{\eta}^{A^{\alpha}h}(r,s)|^2dx\bigr| = \bigl|4\int_{\mathcal{D}}P_N\mathcal{O}_r^{\varrho}\tilde{\eta}^{A^{\alpha}h}(r,s)\partial_x\tilde{\eta}^{A^{\alpha}h}(r,s)dx\bigr|\\
		\le&\ 4\|P_N\mathcal{O}_r^{\varrho}\|_V\|\tilde{\eta}^{A^{\alpha}h}(r,s)\|\cdot\|\tilde{\eta}^{A^{\alpha}h}(r,s)\|_1 \le  \frac{2}{\rho}\|P_N\mathcal{O}_r^{\varrho}\|^2_V\|\tilde{\eta}^{A^{\alpha}h}(r,s)\|^2+2\rho\|\tilde{\eta}^{A^{\alpha}h}(r,s)\|_1^2.
	\end{align*}
	By the integration by parts formula, the H\"older inequality, the Gagliardo-Nirenberg inequality, and $\rho$-Young's inequalities $ab\le\frac{3}{4(4\rho)^{\frac{1}{3}}}a^{\frac{4}{3}}+\rho b^4$ and $\frac{3}{4\rho^{\frac{1}{3}}}a^{\frac{4}{3}} \le \frac{\varepsilon}{2}a^2+\frac{1}{4\rho\varepsilon^2}$, we get
	\begin{align*}
		III=&\ \bigl|2\int_{\mathcal{D}}Y^{N,\varrho}(r)\partial_x\tilde{\eta}^{A^{\alpha}h}(r,s)\cdot \tilde{\eta}^{A^{\alpha}h}(r,s)dx\bigr|\\ 
		=&\ \bigl|\int_{\mathcal{D}}Y^{N,\varrho}(r)\partial_x\big(\tilde{\eta}^{A^{\alpha}h}(r,s)\big)^2dx\bigr| = \bigl|\int_{\mathcal{D}}\partial_xY^{N,\varrho}(r)|\tilde{\eta}^{A^{\alpha}h}(r,s)|^2dx\bigr|\\
		\le&\ \sqrt{2}\|\nabla Y^{N,\varrho}(r)\|\cdot\|\tilde{\eta}^{A^{\alpha}h}(r,s)\|^{\frac{3}{2}}\|\tilde{\eta}^{A^{\alpha}h}(r,s)\|_1^{\frac{1}{2}}\\
		\le&\ \frac{3}{4\rho^{\frac{1}{3}}}\|\nabla Y^{N,\varrho}(r)\|^{\frac{4}{3}}\|\tilde{\eta}^{A^{\alpha}h}(r,s)\|^2 + \rho\|\tilde{\eta}^{A^{\alpha}h}(r,s)\|^2_1\\
		\le&\ \bigl(\frac{\varepsilon}{2}\|\nabla Y^{N,\varrho}(r)\|^{2}+\frac{1}{4\rho\varepsilon^2}\bigr)\|\tilde{\eta}^{A^{\alpha}h}(r,s)\|^2 + \rho\|\tilde{\eta}^{A^{\alpha}h}(r,s)\|^2_1,
	\end{align*}
	and 
	\begin{align*}
		IV=&\ \bigl|2\int_{\mathcal{D}}P_N\mathcal{O}_r^{\varrho}\partial_x\tilde{\eta}^{A^{\alpha}h}(r,s)\cdot \tilde{\eta}^{A^{\alpha}h}(r,s)dx\bigr| \le 2\|P_N\mathcal{O}_r^{\varrho}\|_V\|\partial_x\tilde{\eta}^{A^{\alpha}h}(r,s)\|\|\tilde{\eta}^{A^{\alpha}h}(r,s)\|\\
		\le&\ \rho\|\tilde{\eta}^{A^{\alpha}h}(r,s)\|_1^2 + \frac{1}{\rho}\|P_N\mathcal{O}_r^{\varrho}\|_V^2\|\tilde{\eta}^{A^{\alpha}h}(r,s)\|^2.
	\end{align*}
	By Lemma \ref{l11.2.1}, we have
	\begin{align*}
		&\ \bigl|\bigl(P_NDF(X^N(r))S(r-s)A^{\alpha}h,\tilde{\eta}^{A^{\alpha}h}(r,s)\bigr)\bigr|\\
		=&\ \bigl|\bigl(A^{-\frac{1}{2}}DF(X^N(r))S(r-s)A^{\alpha}h,A^{\frac{1}{2}}\tilde{\eta}^{A^{\alpha}h}(r,s)\bigr)\bigr|\\
		\le&\ \frac{1}{2}\|DF(X^N(r))S(r-s)A^{\alpha}h\|_{-1}^2 + \frac{1}{2}\|\tilde{\eta}^{A^{\alpha}h}(r,s)\|_1^{2}\\
		\le&\ C\|X^N(r)\|_V^2\|S(r-s)A^{\alpha}h\|^2 + \frac{1}{2}\|\tilde{\eta}^{A^{\alpha}h}(r,s)\|_1^{2}.
	\end{align*}	
	Above all, we obtain
	\begin{align*}
		&\ \frac{1}{2}\|\tilde{\eta}^{A^{\alpha}h}(t,s)\|^2 + \int_s^t\|\tilde{\eta}^{A^{\alpha}h}(r,s)\|_1^2dr \le \bigl(\frac{3}{4}+4\rho\bigr)\int_s^t\|\tilde{\eta}^{A^{\alpha}h}(r,s)\|_1^2dr\\ 
		&\ + \int_s^t\bigl(\varepsilon\|\nabla Y^{N,\varrho}(r)\|^{2} + \frac{16}{\varepsilon^2} + \frac{1}{4\rho\varepsilon^2} + \frac{3}{\rho}\|P_N\mathcal{O}_r^{\varrho}\|^2_V\bigr)\|\tilde{\eta}^{A^{\alpha}h}(r,s)\|^2dr\\ 
		&\ + C\int_s^t\|X^N(r)\|_V^2\|S(r-s)A^{\alpha}h\|^2dr.
	\end{align*}
	Note that $\frac{3}{4}+4\rho<1$ and by \eqref{eq3.2.2}, then we get
	\begin{align*}
		\|\tilde{\eta}^{A^{\alpha}h}(t,s)\|^2
		\le&\ \int_s^t\big(2\varepsilon\|\nabla Y^{N,\varrho}(r)\|^{2} + \frac{32}{\varepsilon^2} + \frac{1}{2\rho\varepsilon^2} + \frac{6}{\rho}\|P_N\mathcal{O}_r^{\varrho}\|^2_V\big)\|\tilde{\eta}^{A^{\alpha}h}(r,s)\|^2dr\\ 
		&\ + C\sup_{u\in[0,T]}\|X^N(u)\|_V^2\int_s^t\|S(r-s)A^{\alpha}h\|^2dr\\
		\le&\ \int_s^t\big(2\varepsilon\|\nabla Y^{N,\varrho}(r)\|^{2} + C(\rho,\varepsilon) + \frac{6}{\rho}\|P_N\mathcal{O}_r^{\varrho}\|^2_V\big)\|\tilde{\eta}^{A^{\alpha}h}(r,s)\|^2dr\\ 
		&\ + C\sup_{u\in[0,T]}\|X^N(u)\|_V^2\cdot C(t-s)^{1-2\alpha}\|h\|^2.
	\end{align*}
	From a Gronwall's inequality, it follows that
	\begin{align*}
		&\ \|\tilde{\eta}^{A^{\alpha}h}(t,s)\|^2\\ \le&\  CT^{1-2\alpha}\sup_{u\in[0,T]}\|X^N(u)\|_V^2\|h\|^2\exp\Bigl(\int_s^t2\varepsilon\|\nabla Y^{N,\varrho}(r)\|^{2} + C(\rho,\varepsilon) + \frac{6}{\rho}\|P_N\mathcal{O}_r^{\varrho}\|^2_Vdr\Bigr)\\
		\le&\ C\sup_{u\in[0,T]}\|X^N(u)\|_V^2\|h\|^2\exp\Bigl(\int_s^t2\varepsilon\|\nabla Y^{N,\varrho}(r)\|^{2} + \frac{6}{\rho}\|P_N\mathcal{O}_r^{\varrho}\|^2_Vdr\Bigr).
	\end{align*}
	Thus we have
	\begin{equation}\label{eq11.3.7}
		\begin{aligned}
			&\ \|\Pi(t,s)A^{\alpha}h\| = \|\eta^{A^{\alpha}h}(t,s)\| \le \|\tilde{\eta}^{A^{\alpha}h}(t,s)\| + \|S(t-s)A^{\alpha}h\|\\ 
			\le&\ C\sup_{u\in[0,T]}\|X^N(u)\|_V\|h\|\exp\Bigl(\int_s^t\varepsilon\|\nabla Y^{N,\varrho}(r)\|^{2} + \frac{3}{\rho}\|P_N\mathcal{O}_r^{\varrho}\|^2_Vdr\Bigr) + C(t-s)^{-\alpha}\|h\|.
		\end{aligned}
	\end{equation}
	
	As $\alpha\in(\frac{1}{2},1)$, according to the Duhamel formula, the solution of \eqref{eq14.4.3} is of the form
	\begin{equation}\label{eq11.3.6}
		\tilde{\eta}^{A^{\alpha}h}(t,s) = \int_s^t\Pi(t,r)P_NDF(X^N(r))S(r-s)A^{\alpha}hdr.
	\end{equation}
	By \eqref{eq11.3.6}, \eqref{eq11.3.7}, Lemma \ref{l11.2.1}, and Sobolev embedding $\dot{H}^1\hookrightarrow V$, we have
	\begin{align*}
		&\|\tilde{\eta}^{A^{\alpha}h}(t,s)\| \le  \int_s^t\!\|\Pi(t,r)A^{\frac{1}{2}}\big[A^{-\frac{1}{2}}DF(X^N(r))S(r\!-\!s)A^{\alpha}h\big]\|dr\\
		\le&  \int_s^tC\sup_{u\in[0,T]}\|X^N(u)\|_V\|A^{-\frac{1}{2}}P_NDF(X^N(r))S(r-s)A^{\alpha}h\| \exp\Bigl(\int_r^t \varepsilon\|\nabla Y^{N,\varrho}(u)\|^{2} + \frac{3}{\rho}\|P_N\mathcal{O}_u^{\varrho}\|^2_Vdu\Bigr)dr\\
		& + \int_s^tC(t-r)^{-\frac{1}{2}}\|A^{-\frac{1}{2}}DF(X^N(r))S(r-s)A^{\alpha}h\|dr\\
		\le& C\sup_{u\in[0,T]}\|X^N(u)\|_V\int_s^tC\|X^N(r)\|_V(r-s)^{-\alpha}\|h\|\exp\Bigl(\int_r^t \varepsilon\|\nabla Y^{N,\varrho}(u)\|^{2} + \frac{3}{\rho}\|P_N\mathcal{O}_u^{\varrho}\|^2_Vdu\Bigr)dr\\
		& +\int_s^tC(t-r)^{-\frac{1}{2}}\|X^N(r)\|_V(r-s)^{-\alpha}\|h\|dr\\
		\le& C\sup_{u\in[0,T]}\|X^N(u)\|^2_V\int_s^t(r-s)^{-\alpha}\|h\|\exp\Bigl(\int_r^t \varepsilon\|\nabla Y^{N,\varrho}(u)\|^{2} + \frac{3}{\rho}\|P_N\mathcal{O}_u^{\varrho}\|^2_Vdu\Bigr)dr\\ 
		&\ +C(t-s)^{\frac{1}{2}-\alpha}\sup_{u\in[0,T]}\|X^N(u)\|_V\|h\|,
	\end{align*}
	where we have used the fact that $\int_s^t(t-r)^{-\frac12}(r-s)^{-\alpha}\,dr = \mathcal B\left(\frac12,1-\alpha\right)(t-s)^{\frac12-\alpha}$ for all $0\le\alpha<1$ and $0\le s<t\le T$.
	
	Therefore, by \eqref{eq3.2.2} we obtain
	\begin{align*}
		&\ \|\Pi(t,s)A^{\alpha}h\| 
		= \|\eta^{A^{\alpha}h}(t,s)\| \le \|\tilde{\eta}^{A^{\alpha}h}(t,s)\| 
		+ \|S(t-s)A^{\alpha}h\| \\
		\le&\ C\sup_{u\in[0,T]}\|X^N(u)\|^2_V\int_s^t(r-s)^{-\alpha}\|h\|\exp\Bigl(\int_r^t \varepsilon\|\nabla Y^{N,\varrho}(u)\|^{2} + \frac{3}{\rho}\|P_N\mathcal{O}_u^{\varrho}\|^2_Vdu\Bigr)dr\\ 
		&\ +C(t-s)^{-\alpha}\sup_{u\in[0,T]}\|X^N(u)\|_V\|h\| + C(t-s)^{-\alpha}\|h\|.
	\end{align*}
	This completes the proof.
\end{proof}

Now we derive moments estimates for $\eta^h(t)$.

\begin{lemma}\label{l14.3.1.1}
	Assume (A1) and (A2). Then for any $\alpha\in[0,1)$, $\delta\in[0,\frac{1}{2})$ and $q\ge1$, there exists a constant $C > 0$ independent of $N$ such that for all $h \in H_N$ and $0\le t\le T$,
	\begin{align*}
		\|\eta^{A^{\alpha}h}(t)\|_{L^{2q}(\Omega;H)} \le&\ C\bigl(1+t^{-\alpha}+\|P_N\xi\|^4_{V}\bigr)e^{\frac{\lambda_0}{4}\|\xi\|^2}\|h\|,\\
		\|A^{\delta}\eta^{A^{\alpha}h}(t)\|_{L^{2q}(\Omega;H)} \le&\ C(1+t^{-(\delta+\alpha)}+\|P_N\xi\|_{V}^6)e^{\frac{\lambda_0}{4}\|\xi\|^2}\|h\|.
	\end{align*}
\end{lemma}

\begin{proof}
	Take $0<\varepsilon'<1$ and $G(\omega) = \big(K_{\varepsilon,p_0}^{(\kappa)}(\omega)\big)^{\big(\frac{1}{4}-\varepsilon-\gamma-\kappa\big)^{-1}}>0$ with $\varepsilon+\kappa<\frac{1}{4}-\gamma$, then by Lemma \ref{l14.3.2} we have
	\begin{equation*}
		\varrho_0(\omega) = \frac{1}{\varepsilon'} \Bigl(C_T\big(K_{\varepsilon_1,4}(\omega)\big)^{(\frac{1}{4}-\varepsilon_1-\gamma)^{-1}}+1+ \big(K_{\varepsilon,p_0}^{(\kappa)}(\omega)\big)^{\big(\frac{1}{4}-\varepsilon-\gamma-\kappa\big)^{-1}}\Bigr) > 1.
	\end{equation*}
	Then
	\begin{equation*}
		\varepsilon'\varrho_0(\omega) \ge \big(K_{\varepsilon,p_0}^{(\kappa)}(\omega)\big)^{\big(\frac{1}{4}-\varepsilon-\gamma-\kappa\big)^{-1}},
	\end{equation*}
	which implies
	\begin{equation*}
		K_{\varepsilon,p_0}^{(\kappa)}(\omega) \le \big(\varepsilon'\varrho_0(\omega)\big)^{\frac{1}{4}-\varepsilon-\gamma-\kappa}.
	\end{equation*}
	Therefore, by \eqref{eq14.2.9.1} and $\|P_Nu\|_{L^{p_0}(\mathcal{D})} \le C\|u\|_{L^{p_0}(\mathcal{D})}$ for all $u\in L^{p_0}(\mathcal{D})$, we obtain
	\begin{align*}
		&\ \|P_N\mathcal{O}_u^{\varrho_0(\omega)}\|_V \le C\|A^{\kappa}P_N\mathcal{O}_u^{\varrho_0(\omega)}\|_{L^{p_0}(\mathcal{D})} \le C\|A^{\kappa}\mathcal{O}_u^{\varrho_0(\omega)}\|_{L^{p_0}(\mathcal{D})}\\ 
		\le&\ C\cdot C_T\varrho_0^{-\frac{1}{4}+\varepsilon+\gamma+\kappa}K_{\varepsilon,p_0}^{(\kappa)} \le C(\varepsilon')^{\frac{1}{4}-\varepsilon-\gamma-\kappa} \le C_T.
	\end{align*}
	By Lemma \ref{t14.3.1} and taking $0<\varepsilon\le \frac{\lambda_0}{16q}$, we get
	\begin{equation}\label{eq14.3.6}
		\begin{aligned}
			&\ \mathbb{E}\Bigl[\exp\Bigl(\int_r^t 4q\varepsilon\|\nabla Y^{N,\varrho}(u)\|^{2} + \frac{12q}{\rho}\|P_N\mathcal{O}_u^{\varrho}\|^2_Vdu\Bigr)\Bigr]\\
			\le&\ \mathbb{E}\Bigl[\exp\Bigl(\int_r^t \frac{\lambda_0}{4}\|\nabla Y^{N,\varrho}(u)\|^{2}du\Bigr)\cdot\exp\Bigl(\int_r^t\frac{12q}{\rho}\|P_N\mathcal{O}_u^{\varrho}\|^2_Vdu\Bigr)\Bigr]\\
			\le&\ \mathbb{E}\Bigl[\exp\Bigl(\int_r^t \frac{\lambda_0}{4}\|\nabla Y^{N,\varrho}(u)\|^{2}du\Bigr)\cdot\exp\Bigl(\frac{12q}{\rho}C_T^2T\Bigr)\Bigr] \le C\exp\bigl[\frac{\lambda_0}{4}\|\xi\|^2\bigr].
		\end{aligned}
	\end{equation}
	Note that $\eta^{A^{\alpha}h}(t) = \Pi(t,0)A^{\alpha}h$, by the H\"older inequality and \eqref{eq14.3.6}, we have
	\begin{align*}
		&\ \|\eta^{A^{\alpha}h}(t)\|_{L^{2q}(\Omega;H)} = \bigl(\mathbb{E}\|\eta^{A^{\alpha}h}(t)\|^{2q}\bigr)^{\frac{1}{2q}}\\ 
		\le&\ C\Bigl(\mathbb{E}\Bigl[\sup_{u\in[0,T]}\|X^N(u)\|^2_V\int_0^tr^{-\alpha}\|h\|\exp\Bigl(\int_r^t \varepsilon\|\nabla Y^{N,\varrho}(u)\|^{2} + \frac{3}{\rho}\|P_N\mathcal{O}_u^{\varrho}\|^2_Vdu\Bigr)dr\Bigr]^{2q}\Bigr)^{\frac{1}{2q}}\\ 
		&\ +Ct^{-\alpha}\Bigl(\mathbb{E}\Bigl[\sup_{u\in[0,T]}\|X^N(u)\|_V\|h\|\Bigr]^{2q}\Bigr)^{\frac{1}{2q}} + Ct^{-\alpha}\|h\|\\
		\le&\ C\bigl[\mathbb{E}\sup_{u\in[0,T]}\|X^N(u)\|_V^{8q}\bigr]^{\frac{1}{4q}}\int_0^tr^{-\alpha}\|h\|\\
		&\ \cdot \Big(\mathbb{E}\Big[\exp\Big(\int_r^t4q\varepsilon\|\nabla Y^{N,\varrho}(u)\|^{2} + \frac{12q}{\rho}\|P_N\mathcal{O}_u^{\varrho}\|^2_Vdu\Big)\Big]\Big)^{\frac{1}{4q}}dr + Ct^{-\alpha}(1+\|P_N\xi\|_V^2)\|h\| + Ct^{-\alpha}\|h\|\\
		\le&\ C(1+\|P_N\xi\|^4_{V})\int_0^tr^{-\alpha}\|h\|e^{\frac{\lambda_0}{4}\|\xi\|^2}dr +Ct^{-\alpha}(1+\|P_N\xi\|_V^2)\|h\| + Ct^{-\alpha}\|h\|\\
		\le&\  C(1+t^{-\alpha}+\|P_N\xi\|^4_{V})e^{\frac{\lambda_0}{4}\|\xi\|^2}\|h\|.
	\end{align*}
	This completes the first assertion.
	
	Note that the mild solution of \eqref{eq14.4.2} with initial value $A^{\alpha}h$ is of the form
	\begin{equation*}
		\eta^{A^{\alpha}h}(t) = S(t)A^{\alpha}h + \int_0^tS(t-r)P_NDF(X^N(r))\eta^{A^{\alpha}h}(r)dr.
	\end{equation*}
	For any $\delta \in[0,\frac{1}{2})$, by \eqref{eq3.2.2} and Lemma \ref{l11.2.1}, we have
	\begin{align*}
		\|A^{\delta}\eta^{A^{\alpha}h}(t)\|
		\le&\ \|A^{\delta}S(t)A^{\alpha}h\| + \int_0^t\|A^{\delta+\frac{1}{2}}S(t-r)\bigl[A^{-\frac{1}{2}}P_NDF(X^N(r))\eta^{A^{\alpha}h}(r)\bigr]\|dr\\
		\le&\ Ct^{-(\delta+\alpha)}\|h\| + \int_0^tC(t-r)^{-\frac{1}{2}-\delta}C\|A^{-\frac{1}{2}}P_NDF(X^N(r))\eta^{A^{\alpha}h}(r)\|dr\\
		\le&\ Ct^{-(\delta+\alpha)}\|h\| + \int_0^tC(t-r)^{-\frac{1}{2}-\delta}\|X^N(r)\|_V\|\eta^{A^{\alpha}h}(r)\|dr.
	\end{align*}
	By the H\"older inequality and Theorem \ref{t14.2.6}, we obtain
	\begin{align*}
		&\ \|A^{\delta}\eta^{A^{\alpha}h}(t)\|_{L^{2q}(\Omega;H)}\\
		\le&\ Ct^{-(\delta+\alpha)}\|h\| + \int_0^tC(t-r)^{-\delta-\frac{1}{2}}\|X^N(r)\|_{L^{4q}(\Omega;V)}\|\eta^{A^{\alpha}h}(r)\|_{L^{4q}(\Omega;H)}dr\\
		\le&\ Ct^{-(\delta+\alpha)}\|h\| + \int_0^tC(t-r)^{-\delta-\frac{1}{2}}\cdot C(1+\|P_N\xi\|^2_{V})\cdot C(1+r^{-\alpha}+\|P_N\xi\|^4_{V})e^{\frac{\lambda_0}{4}\|\xi\|^2}\|h\|dr\\
		\le&\ C(1+t^{-(\delta+\alpha)}+\|P_N\xi\|_{V}^6)e^{\frac{\lambda_0}{4}\|\xi\|^2}\|h\|.
	\end{align*}
	Thus the proof is completed.
\end{proof}

Now we investigate estimates for the first order Fr\'echet  derivative $DU^N(t,\nu)$.

\begin{theorem}\label{t12.3.1}
	Assume (A1) and (A2). Then for any $\alpha\in[0,1)$, there exists a constant $C>0$ independent of $N$ such that
	\begin{equation*}
		\|A^{\alpha}DU^N(t,\xi)\| \le C(1+t^{-\alpha}+\|P_N\xi\|^4_{V})e^{\frac{\lambda_0}{4}\|\xi\|^2}.
	\end{equation*}
\end{theorem}

\begin{proof}
	For any $\alpha\in[0,1)$ and $h\in H$, by Lemma \ref{l14.4.2} we have
	\begin{align*}
		|(A^{\alpha}DU^N(t,\xi),h)|=&\ |\mathbb{E}(D\phi(X^N(t,\xi))DX^N(t,\xi),A^{\alpha}h)| = |\mathbb{E}[D\phi(X^N(t,\xi))(DX^N(t,\xi)A^{\alpha}h)]|\\
		=&\ |\mathbb{E}(D\phi(X^N(t,\xi)),\eta^{A^{\alpha}h}(t))| \le \mathbb{E}\big[\|D\phi(X^N(t,\xi))\|\cdot\|\eta^{A^{\alpha}h}(t)\|\big]\\
		\le&\ |\phi|_{\mathcal{C}_b^1(H)}\|\eta^{A^{\alpha}h}(t)\|_{L^1(\Omega;H)} \le C(1+t^{-\alpha}+\|P_N\xi\|^4_{V})e^{\frac{\lambda_0}{4}\|\xi\|^2}\|h\|.
	\end{align*}
	Thus the proof is completed.
\end{proof}

For any $h,g\in H_N$, the second order directional Fr\'echet derivative $D^2U^N(t,\xi)(h,g)$ of $U^N(t,\xi)$ with respect to $\xi$ along directions $h$ and $g$ is of the form, cf.\ \cite{ZhangZou2026JSC},
\begin{equation*}
	D^2U^N(t,\xi)(h,g) = \mathbb{E}\big[D^2\phi(X^N(t,\xi))\big(\eta^h(t),\eta^g(t)\big)\big] + \mathbb{E}\big[D\phi(X^N(t,\xi))\zeta^{h,g}(t)\big],
\end{equation*}
where $\zeta^{h,g}(t) = D^2X^N(t,\xi)(h,g)$ satisfies $\zeta^{h,g}(0)=0$ and
\begin{equation}\label{eq11.3.4}
	d\zeta^{h,g}(t) + A\zeta^{h,g}(t)dt = P_NDF(X^N(t))\zeta^{h,g}(t)dt + P_ND^2F(X^N(t))(\eta^h(t),\eta^g(t))dt.
\end{equation}
The next lemma studies an estimate for $\zeta^{h,g}(t)$.

\begin{lemma}\label{l12.3.2}
	Assume (A1) and (A2). Then for any $\alpha_1,\alpha_2\ge0$ with $\alpha_1+\alpha_2\in[0,1)$ and $q\ge1$, there exists a constant $C > 0$ independent of $N$ such that for all $h,g \in H$ and $t \in [0,T]$,
	\begin{equation*}
		\|\zeta^{A^{\alpha_1}h,A^{\alpha_2}g}(t)\|_{L^{2q}(\Omega;H)} \le C(1+t^{-\alpha_1-\alpha_2}+\|P_N\xi\|_{V}^{16})e^{\frac{3}{4}\lambda_0\|\xi\|^2}\|h\|\|g\|.
	\end{equation*}
\end{lemma}

\begin{proof}
	According to the Duhamel formula, the solution to \eqref{eq11.3.4} is formulated as
	\begin{equation*}
		\zeta^{h,g}(t) = \int_0^t\Pi(t,s)P_ND^2F(X^N(s))\big(\eta^h(s),\eta^g(s)\big)ds.
	\end{equation*}
	Then, we have
	\begin{equation*}
		\zeta^{A^{\alpha_1}h,A^{\alpha_2}g}(t) = \int_0^t\Pi(t,s)P_ND^2F(X^N(s))\big(\eta^{A^{\alpha_1}h}(s),\eta^{A^{\alpha_2}g}(s)\big)ds.
	\end{equation*}
	To simplify notation, define $\mathcal{B} = A^{-\frac{3}{4}-\rho}P_ND^2F(X^N(s))\bigl(\eta^{A^{\alpha_1}h}(s),\eta^{A^{\alpha_2}g}(s)\bigr)$ for any  $0<\rho<\min\{\frac{1}{16},\frac{1-\alpha_1-\alpha_2}{2}\}$. By Sobolev embedding $L^1(\mathcal{D})\hookrightarrow \dot{H}^{-\frac{1}{2}-2\rho}$, Lemma \ref{l11.2.1}, \eqref{eq12.2.3}, and H\"older inequality, we get for $p=\frac{1}{1-4\rho}>1$,
	\begin{align*}
		\|\mathcal{B}\|_{L^{8q}(\Omega;H)} =&\ \|A^{-\frac{1}{4}-\rho}A^{-\frac{1}{2}}P_ND^2F(X^N(s))\bigl(\eta^{A^{\alpha_1}h}(s),\eta^{A^{\alpha_2}g}(s)\bigr)\|_{L^{8q}(\Omega;H)}\\
		\le&\ \|A^{-\frac{1}{2}}D^2F(X^N(s))\bigl(\eta^{A^{\alpha_1}h}(s),\eta^{A^{\alpha_2}g}(s)\bigr)\|_{L^{8q}(\Omega;L^1(\mathcal{D}))}\\
		\le&\ C\bigl\|\|\eta^{A^{\alpha_1}h}(s)\|_{L^{2p}(\mathcal{D})}\|\eta^{A^{\alpha_2}g}(s)\|_{L^{2p}(\mathcal{D})}\bigr\|_{L^{8q}(\Omega)}\\
		\le&\ C\|\eta^{A^{\alpha_1}h}(s)\|_{L^{16q}(\Omega;L^{2p}(\mathcal{D}))}\|\eta^{A^{\alpha_2}g}(s)\|_{L^{16q}(\Omega;L^{2p}(\mathcal{D}))}.
	\end{align*}
	Note that $\alpha_1+\alpha_2+\frac{1}{2}-\frac{1}{2p}=\alpha_1+\alpha_2+2\rho<1$. By Sobolev embedding inequality $\dot{H}^{2\rho}=\dot{H}^{\frac{1}{2}-\frac{1}{2p}}\hookrightarrow L^{2p}(\mathcal{D})$ and Lemma \ref{l14.3.1.1}, we get
	\begin{align*}
		\|\mathcal{B}\|_{L^{8q}(\Omega;H)} \le&\ C\|A^{\rho}\eta^{A^{\alpha_1}h}(s)\|_{L^{16q}(\Omega;H)}\|A^{\rho}\eta^{A^{\alpha_2}g}(s)\|_{L^{16q}(\Omega;H)}\\
		\le&\ C(1+s^{-(\rho+\alpha_1)}+\|P_N\xi\|_{V}^6)e^{\frac{\lambda_0}{4}\|\xi\|^2}\|h\| \cdot C(1+s^{-(\rho+\alpha_2)}+\|P_N\xi\|_{V}^6)e^{\frac{\lambda_0}{4}\|\xi\|^2}\|g\|\\
		\le&\ C(1+s^{-2\rho-\alpha_1-\alpha_2})(1+\|P_N\xi\|_{V}^{12})e^{\frac{1}{2}\lambda_0\|\xi\|^2}\|h\|\|g\|.
	\end{align*}
	By Lemma \ref{l14.4.2},\eqref{eq14.3.6}, and the H\"older inequality, we obtain
	\begin{align*}
		&\ \|\zeta^{A^{\alpha_1}h,A^{\alpha_2}g}(t)\|_{L^{2q}(\Omega;H)} \le \int_0^t\|\Pi(t,s)A^{\frac{3}{4}+\rho}\mathcal{B}\|_{L^{2q}(\Omega;H)}ds\\
		\le&\ \int_0^t\Bigl\|C\sup_{u\in[0,T]}\|X^N(u)\|_V^2\int_s^t(r-s)^{-\frac{3}{4}-\rho}\|\mathcal{B}\|\exp\Bigl(\int_r^t\Bigl[\varepsilon\|\nabla Y^{N,\varrho}(u)\|^2+\frac{3}{\rho}\|P_N\mathcal{O}_u^\varrho\|_V^2\Bigr]du\Bigr)dr\\
		&\ +C(t-s)^{-\frac{3}{4}-\rho}\sup_{u\in[0,T]}\|X^N(u)\|_V\|\mathcal{B}\|+C(t-s)^{-\frac{3}{4}-\rho}\|\mathcal{B}\|\Bigr\|_{L^{2q}(\Omega)}ds\\
		\le&\ C\int_0^t\Bigl\|\sup_{u\in[0,T]}\|X^N(u)\|_V^2\Bigr\|_{L^{8q}(\Omega)}\int_s^t(r-s)^{-\frac{3}{4}-\rho}\|\mathcal{B}\|_{L^{8q}(\Omega;H)}\\
		&\ \quad\quad\cdot\Bigl[\mathbb{E}\exp\Bigl(\int_r^t\Bigl[4q\varepsilon\|\nabla Y^{N,\varrho}(u)\|^2+\frac{12q}{\rho}\|P_N\mathcal{O}_u^\varrho\|_V^2\Bigr]du\Bigr)\Bigr]^{\frac{1}{4q}}drds\\
		&\ +C\int_0^t(t-s)^{-\frac{3}{4}-\rho}\Bigl[\Bigl\|\sup_{u\in[0,T]}\|X^N(u)\|_V\Bigr\|_{L^{4q}(\Omega)}\|\mathcal{B}\|_{L^{4q}(\Omega;H)}+\|\mathcal{B}\|_{L^{2q}(\Omega;H)}\Bigr]ds\\
		\le&\ C(1+\|P_N\xi\|_V^2)^2(1+\|P_N\xi\|_V^{12})e^{(\frac{1}{2}+\frac{1}{16q})\lambda_0\|\xi\|^2}\|h\|\|g\|\int_0^t\int_s^t(r-s)^{-\frac{3}{4}-\rho}(1+s^{-2\rho-\alpha_1-\alpha_2})drds\\
		&\ +C(1+\|P_N\xi\|_V^2)(1+\|P_N\xi\|_V^{12})e^{\frac{1}{2}\lambda_0\|\xi\|^2}\|h\|\|g\|\int_0^t(t-s)^{-\frac{3}{4}-\rho}(1+s^{-2\rho-\alpha_1-\alpha_2})ds.
	\end{align*}
	Since $0<\rho<\min\{\frac{1}{16},\frac{1-\alpha_1-\alpha_2}{2}\}$, we have
	\begin{align*}
		&\ \int_0^t\int_s^t(r-s)^{-\frac{3}{4}-\rho}(1+s^{-2\rho-\alpha_1-\alpha_2})drds\\
		=&\ \frac{t^{\frac{5}{4}-\rho}}{(\frac{1}{4}-\rho)(\frac{5}{4}-\rho)}+\frac{\mathcal{B}(1-2\rho-\alpha_1-\alpha_2,\frac{5}{4}-\rho)}{\frac{1}{4}-\rho}t^{\frac{5}{4}-3\rho-\alpha_1-\alpha_2}\\
		\le&\ C(1+t^{-\alpha_1-\alpha_2}),
	\end{align*}
	and
	\begin{align*}
		&\ \int_0^t(t-s)^{-\frac{3}{4}-\rho}(1+s^{-2\rho-\alpha_1-\alpha_2})ds\\
		=&\ \frac{t^{\frac{1}{4}-\rho}}{\frac{1}{4}-\rho}+\mathcal{B}(1-2\rho-\alpha_1-\alpha_2,\frac{1}{4}-\rho)t^{\frac{1}{4}-3\rho-\alpha_1-\alpha_2}\\
		\le&\ C(1+t^{-\alpha_1-\alpha_2}),\qquad 0<t\le T.
	\end{align*}
	Therefore, noting that $q\ge1$, we obtain
	\begin{align*}
		\|\zeta^{A^{\alpha_1}h,A^{\alpha_2}g}(t)\|_{L^{2q}(\Omega;H)}\le C(1+t^{-\alpha_1-\alpha_2}+\|P_N\xi\|_{V}^{16})e^{\frac{3}{4}\lambda_0\|\xi\|^2}\|h\|\|g\|,\qquad 0<t\le T.
	\end{align*}
	Thus the proof is completed.
\end{proof}

The next Theorem studies the boundedness and Sobolev regularity of the seconded Fr\'echet derivative $D^2U(t,\xi)$.

\begin{theorem}\label{t12.3.2}
	Assume (A1) and (A2). Then for any $\alpha_1,\alpha_2\ge0$ with $\alpha_1+\alpha_2\in[0,1)$, there exists a constant $C>0$ independent of $N$ such that
	\begin{equation*}
		\|A^{\alpha_1}D^2U^N(t,\xi)A^{\alpha_2}\|_{\mathcal{L}(H)} \le C(1+t^{-\alpha_1-\alpha_2}+\|P_N\xi\|_{V}^{16})e^{\frac{3}{4}\lambda_0\|\xi\|^2}.
	\end{equation*}
\end{theorem}

\begin{proof}
	By Lemma \ref{l12.3.2}, we obtain
	\begin{align*}
		&\ |(A^{\alpha_1}D^2U^N(t,\xi)A^{\alpha_2}h,g)| = |(D^2U^N(t,\xi)A^{\alpha_2}h,A^{\alpha_1}g)|\\
		\le&\ |\mathbb{E}D^2\phi(X^N(t,\xi))\big(\eta^{A^{\alpha_2}h}(t),\eta^{A^{\alpha_1}g}(t)\big)| + |\mathbb{E}[D\phi(X^N(t,\xi))\zeta^{A^{\alpha_2}h,A^{\alpha_1}g}(t)]|\\
		\le&\ |\phi|_{\mathcal{C}_b^2(H)}\mathbb{E}\big[\|\eta^{A^{\alpha_2}h}(t)\|\cdot\|\eta^{A^{\alpha_1}g}(t)\|\big] + |\phi|_{\mathcal{C}_b^1(H)}\mathbb{E}\big[1\cdot\|\zeta^{A^{\alpha_2}h,A^{\alpha_1}g}(t)\|\big]\\
		\le&\ C\|\eta^{A^{\alpha_2}h}(t)\|_{L^2(\Omega;H)}\cdot\|\eta^{A^{\alpha_1}g}(t)\|_{L^2(\Omega;H)} + C\|\zeta^{A^{\alpha_2}h,A^{\alpha_1}g}(t)\|_{L^2(\Omega;H)}\\
		\le&\ C(1+t^{-\alpha_1-\alpha_2}+\|P_N\xi\|_{V}^8)e^{\frac{1}{2}\lambda_0\|\xi\|^2}\|h\|\|g\| + C(1+t^{-\alpha_1-\alpha_2}+\|P_N\xi\|_{V}^{16})e^{\frac{3}{4}\lambda_0\|\xi\|^2}\|h\|\|g\|\\
		\le&\ C(1+t^{-\alpha_1-\alpha_2}+\|P_N\xi\|_{V}^{16})e^{\frac{3}{4}\lambda_0\|\xi\|^2}\|h\|\|g\|.
	\end{align*}
	Thus the proof is completed.
\end{proof}

\section{Weak error estimates}
\label{Sect4}

In this section, we investigate the weak convergence order of the spectral Galerkin approximation.

\begin{theorem}\label{t14.5.2}
	Assume (A1) and (A2). Then for any $\phi \in \mathcal{C}_b^2(H)$ and $0<\epsilon\ll 1-4\gamma$, there exists a constant $C>0$ independent of $N$, such that
	\begin{equation*}
		|\mathbb{E}[\phi(X(T))]-\mathbb{E}[\phi(X^{N}(T))]| \le CN^{-1+4\gamma+\epsilon}.
	\end{equation*}
\end{theorem}

\begin{proof}
	For any $\phi \in \mathcal{C}_b^2(H)$ and $M>N$, we have
	\begin{align*}
		&\ |\mathbb{E}[\phi(X(T))]-\mathbb{E}[\phi(X^{N}(T))]|\\
		\le&\  |\mathbb{E}[\phi(X(T))]-\mathbb{E}[\phi(X^{M}(T))]| + |\mathbb{E}[\phi(X^M(T))]-\mathbb{E}[\phi(X^{N}(T))]|\\
		\le&\ |\phi|_{\mathcal{C}_b^1(H)}(\mathbb{E}\|X(T)-X^M(T)\|^2)^{\frac{1}{2}} + |\mathbb{E}[\phi(X^M(T))]-\mathbb{E}[\phi(X^{N}(T))]|.
	\end{align*}
	For any $\alpha\in[0,1)$, the following estimate will be frequently used in forthcoming analysis
	\begin{equation}\label{eq5.4.1}
		\begin{aligned}
			\|(P_M-P_{N})A^{-\alpha}\|_{\mathcal{L}(H)} \le&\ \|(P_M-I)A^{-\alpha}\|_{\mathcal{L}(H)} + \|(I-P^{N})A^{-\alpha}\|_{\mathcal{L}(H)}\\ 
			\le&\ C(M^{-2\alpha}+N^{-2\alpha}).
		\end{aligned}
	\end{equation}
	Obviously, for any $\mathcal{F}_{T}$-measurable, $H_M$-valued random variable $\eta$, there holds, cf.\ \cite[(5.1)]{Andersson2016}
	\begin{equation*}
		\mathbb{E}[\phi(\eta)] = \mathbb{E}[\mathbb{E}[\phi(\eta)|\mathcal{F}_T]] = \mathbb{E}[\mathbb{E}[\phi(X^M(0,\eta))|\mathcal{F}_T]] = \mathbb{E}[U^M(0,\eta)].
	\end{equation*}
	Noticing that $X^N(t)\in H_N\subset H_M$, we split the weak error as
	\begin{equation}\label{eq3.5.2}
		\begin{aligned}
			&\ |\mathbb{E}[\phi(X^M(T))] - \mathbb{E}[\phi(X^N(T))]| = |\mathbb{E}[U^M(T,P_M\xi)] - \mathbb{E}[U^M(0,X^N(T))]|\\
			\le&\ |\mathbb{E}[U^M(T,P_M\xi)] - \mathbb{E}[U^M(T,P_N\xi)]| + |\mathbb{E}[U^M(T,P_N\xi)] - \mathbb{E}[U^M(0,X^N(T))]|.
		\end{aligned}
	\end{equation}
	
	By mean value theorem, \eqref{eq7.2.4}, \eqref{eq5.4.1}, and Theorem \ref{t12.3.1}, we have
	\begin{align*}
		&\ |\mathbb{E}[U^M(T,P_M\xi)] - \mathbb{E}[U^M(T,P_N\xi)]|\\
		=&\ \big|\mathbb{E}\int_0^1 DU^M(T,\theta P_M\xi+(1-\theta)P_N\xi)(P_M-P_N)\xi d\theta\big|\\
		=&\ \big|\mathbb{E}\int_0^1 \big( DU^M(T,\theta P_M\xi+(1-\theta)P_N\xi),(P_M-P_N)\xi\big) d\theta\big|\\
		=&\ \big|\mathbb{E}\int_0^1 \big(A^{\alpha}DU^M(T,\theta P_M\xi+(1-\theta)P_N\xi),(P_M-P_{N})A^{-\alpha}\xi\big)d\theta\big|\\
		\le&\ \mathbb{E}\int_0^1\|A^{\alpha}DU^M(T,\theta P_M\xi+(1-\theta)P_N\xi)\|\cdot\|(P_M-P_{N})A^{-\alpha}\|_{\mathcal{L}(H)}\|\xi\|d\theta\\
		\le&\ C(1+T^{-\alpha}+\|\theta P_M\xi+(1-\theta)P_N\xi\|_V^4)e^{\frac{\lambda_0}{4}\|\xi\|^2}\cdot C(M^{-2\alpha}+N^{-2\alpha})\|\xi\|\\ 
		\le&\ C(M^{-2\alpha}+N^{-2\alpha}).
	\end{align*}
	
	Now we estimate the second term in \eqref{eq3.5.2}.
	Applying It\^o's formula to $U^M(T-t,X^{N}(t))$ for $t \in [0,T]$, we get
	\begin{align*}
		& U^M(0,X^{N}(T)) - U^M(T,X^{N}(0))\\
		=& \int_0^T -\frac{\partial U^M}{\partial t}(T-t,X^N(t)) + \big(DU^M(T-t,X^N(t)),-AX^N(t)+P_NF(X^N(t))\big)dt\\
		+& \int_0^T\!\big(DU^M(T-t, X^{N}(t)),P_NA^{\gamma}dW(t)\big) + \frac{1}{2}\!\int_0^T\!\mathrm{Tr}\{D^2U^M(T-t,X^{N}(t))P_NA^{2\gamma}P_N\}dt.
	\end{align*}
	Replacing $t$ by $T-t$, $N$ by $M$, and $\nu$ by $X^{N}(t)$ in Kolmogorov equation \eqref{eq14.4.1}, respectively, we get
	\begin{align*}
		\frac{\partial U^M}{\partial t}(T-t,X^{N}(t)) =&\ \big(DU^M(T-t,X^{N}(t)), -AX^N(t)+P_MF(X^{N}(t))\big)\\
		&\ + \frac{1}{2}\mathrm{Tr}\{D^2U^M(T-t,X^{N}(t))P_MA^{2\gamma}P_M\}.
	\end{align*}
	Above two equations imply
	\begin{equation*}
		\begin{aligned}
			&\ |\mathbb{E}[U^M(T,P_N\xi)] - \mathbb{E}[U^M(0,X^N(T))]|\\ 
			\le&\ \big|\mathbb{E}\int_{0}^{T}(DU^M(T-t, X^{N}(t)), (P_M-P_N)F(X^N(t))) dt\big|\\
			&\ + \frac{1}{2}\big| \mathbb{E}\int_0^T \mathrm{Tr}\big\{D^2U^M(T-t,X^{N}(t))\bigl(P_NA^{2\gamma}P_N-P_MA^{2\gamma}P_M\bigr)\big\} dt\big|\\
			=:&\ E^1 + E^2.
		\end{aligned}
	\end{equation*}
	By \eqref{eq5.4.1}, H\"older inequality, Lemma \ref{l14.5.1}, and Theorem \ref{t12.3.1}, we have
	\begin{align*}
		E^1
		\le&\ \mathbb{E}\int_{0}^{T}\big|\big(A^{\alpha}DU^M(T-t, X^{N}(t)), (P_M-P_N)A^{-\alpha+\frac{1}{2}}A^{-\frac{1}{2}}\partial_x(X^N(t)^2)\big)\big| dt\\
		\le&\ \mathbb{E}\int_{0}^{T} \|A^{\alpha}DU^M(T-t, X^{N}(t))\|\cdot\|(P_M-P_N)A^{-\alpha+\frac{1}{2}}\|_{\mathcal{L}(H)}\|A^{-\frac{1}{2}}\partial_x(X^N(t)^2)\|dt\\
		\le&\ \mathbb{E}\int_{0}^{T} C(1+(T-t)^{-\alpha}+\|P_MX^N(t)\|^4_{V})e^{\frac{\lambda_0}{4}\|X^N(t)\|^2}\\
		&\ \quad\quad\cdot C(M^{-2\alpha+1}+N^{-2\alpha+1})\cdot C\|X^N(t)\|_{L^4(\mathcal{D})}^2dt\\
		\le&\ \int_{0}^{T} C\big(\mathbb{E}(1+(T-t)^{-\alpha}+\|X^N(t)\|^{4}_{V})^4\big)^{\frac{1}{4}}\big(\mathbb{E}e^{\lambda_0\|X^N(t)\|^2}\big)^{\frac{1}{4}}\\ 
		&\ \quad\quad\cdot C(M^{-2\alpha+1}+N^{-2\alpha+1}) \big(\mathbb{E}\|X^N(t)\|_{L^4(\mathcal{D})}^4\big)^{\frac{1}{2}}dt\\
		\le&\ C(M^{-2\alpha+1}+N^{-2\alpha+1}).
	\end{align*}
	Notice that
	\begin{align*}
		&\ P_NA^{2\gamma}P_N - P_MA^{2\gamma}P_M = P_NA^{2\gamma}(P_N-P_M) + (P_N-P_M)A^{2\gamma}P_M\\
		=&\ -A^{2\gamma}(P_M-P_N)P_N + (P_N-P_M)A^{2\gamma}P_M
		= (P_N-P_M)A^{2\gamma}P_M.
	\end{align*}
	By \eqref{eq5.4.1} and Theorem \ref{t12.3.1}, we have for small $0<\epsilon\ll 1-4\gamma$,
	\begin{align*}
		E^{2} \le&\ \frac{1}{2} \mathbb{E}\int_0^T\big| \mathrm{Tr}\big\{D^2U^M(T-t,X^{N}(t))(P_N-P_M)A^{2\gamma}P_M\big\}\big| dt\\
		=&\ \frac{1}{2} \mathbb{E}\int_0^T\big| \mathrm{Tr}\big\{A^{\alpha}D^2U^M(T-t,X^{N}(t)) P_M(P_N-P_M)A^{-\alpha+2\gamma+\frac{1}{2}+\epsilon}A^{-\frac{1}{2}-\epsilon}\big\}\big|dt \\
		\le&\ \frac{1}{2} \mathbb{E}\int_0^T\|A^{\alpha}D^2U^M(T-t,X^{N}(t))\|_{\mathcal{L}(H)}\|(P_N-P_M)A^{-\alpha+2\gamma+\frac{1}{2}+\epsilon}\|_{\mathcal{L}(H)}\big| \mathrm{Tr}\big\{A^{-\frac{1}{2}-\epsilon}\big\}\big|dt \\
		\le&\ \frac{1}{2} \mathbb{E}\int_0^T C(1+(T-t)^{-\alpha}+\|P_MX^{N}(t)\|_{V}^{16})e^{\frac{3}{4}\lambda_0\|X^{N}(t)\|^2}\cdot C(M^{-2\alpha+4\gamma+1+2\epsilon}+N^{-2\alpha+4\gamma+1+2\epsilon})\cdot Cdt\\
		\le&\ C(M^{-2\alpha+4\gamma+1+2\epsilon}+N^{-2\alpha+4\gamma+1+2\epsilon}).
	\end{align*}
	Thus we have
	\begin{equation}\label{eq14.4.4}
		\begin{aligned}
			&\ |\mathbb{E}[\phi(X(T))] - \mathbb{E}[\phi(X^N(T))]|\\ 
			\le&\ |\phi|_{\mathcal{C}_b^1(H)}(\mathbb{E}\|X(T)-X^M(T)\|^2)^{\frac{1}{2}} + C(M^{-2\alpha+4\gamma+1+2\epsilon}+N^{-2\alpha+4\gamma+1+2\epsilon}).
		\end{aligned}
	\end{equation}
	
	By Lemma \ref{l14.3.1}, $X^M(T)\to X(T)$ in $H$ $a.s$. Moreover, Theorem \ref{t14.2.6}, Lemma \ref{l14.5.1}, and the uniform bound $\sup_{M\ge1}\|P_M\xi\|_V\le C\|\xi\|_1$ imply that
	\begin{equation*}
		\sup_{M\ge1}\mathbb{E}\|X(T)-X^M(T)\|^4\le8\mathbb{E}\|X(T)\|^4+8\sup_{M\ge1}\mathbb{E}\|X^M(T)\|^4<\infty.
	\end{equation*}
	Hence, the sequence $\{\|X(T)-X^M(T)\|^2\}_{M\ge1}$ is uniformly integrable. By Vitali's convergence theorem, we obtain
	\begin{equation*}
		\lim_{M\to\infty}\mathbb{E}\|X(T)-X^M(T)\|^2=0.
	\end{equation*}
	Given $0<\epsilon<1-4\gamma$, take $\alpha=1-\frac{\epsilon}{4}$ and $\varepsilon=\frac{\epsilon}{4}$. Then $\alpha-2\gamma-\frac{1}{2}-\varepsilon=\frac{1-4\gamma-\epsilon}{2}>0$. Letting $M\to\infty$ in \eqref{eq14.4.4} yields
	\begin{equation*}
		|\mathbb{E}[\phi(X(T))]-\mathbb{E}[\phi(X^N(T))]|\le CN^{-2\alpha+4\gamma+1+2\varepsilon}=CN^{-1+4\gamma+\epsilon}.
	\end{equation*}
	Thus the proof is completed.
\end{proof}

\begin{remark}\label{remark4.3}
	We emphasize that when $\gamma=0$, the noise in \eqref{eq14.1.1} is just the to space-time white noise, and the weak convergence order is arbitrarily close to $1$. On the other hand, $\gamma=\frac14$ represents the critical case, in which the solution to \eqref{eq14.1.1} exists only in the sense of distributions. Consequently, our analysis does not yield any weak convergence result for the numerical approximation in this critical regime.
\end{remark}

\section{Numerical experiment}
\label{Sect5}

\begin{figure}[!h]
	\centering
	\includegraphics[width=0.6\textwidth]{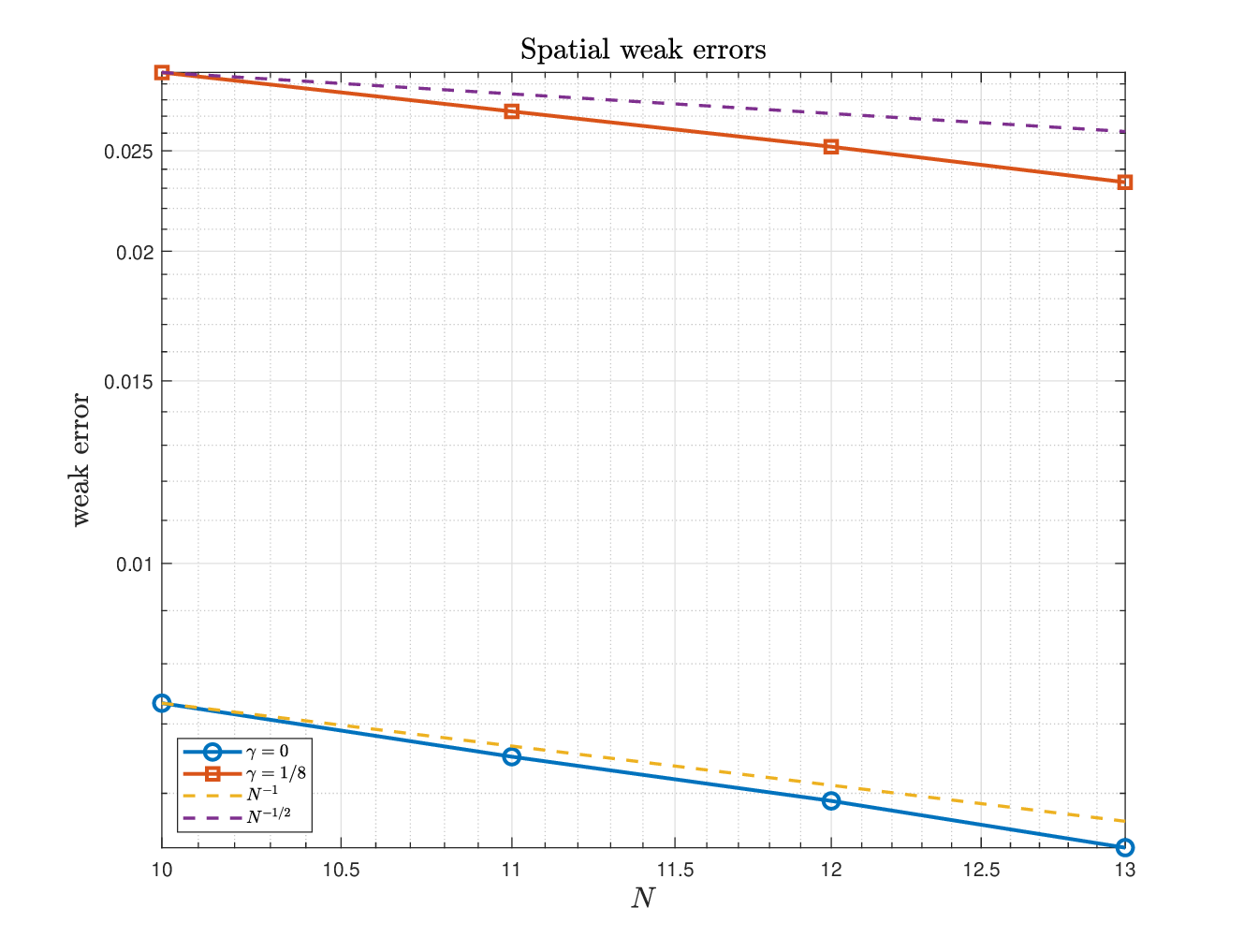}
	\caption{Weak error in log-log scale}
	\label{figure14.1}
\end{figure}

In this section, we carry out a numerical experiment to illustrate theoretical results. 
Consider a stochastic Burgers equation for $(t,x)\in (0,1)\times(0,1)$
\begin{equation*}
	\begin{aligned}
		dX(t,x)-\Delta X(t,x)dt =&\  \nabla(X^2(t,x))dt+(-\Delta)^{\gamma}dW(t),\\
		X(0,x) =&\ \sin(\pi x),
	\end{aligned}
\end{equation*}
where $W(t)$ is a space-time white noise.
We fix time step size $\tau = \frac{1}{1000}$, $N = 10,11,12,13,64$ and apply a tamed exponential Euler method
\begin{equation*}
	X_{k+1}^N = e^{-\tau A}X_k^N + A^{-1}(I-e^{-\tau A})\frac{P_N\partial_x(X_k^N)^2}{1+\tau\|P_N\partial_x(X_k^N)^2\|} + \int_{t_k}^{t_{k+1}}e^{-A(t_{k+1}-s)}P_NA^{\gamma}dW(s),
\end{equation*}
to conduct numerical computations. We choose a test function $\phi(X)=e^{-\|X\|^2}$, Monte Carlo samples $M=2000$, and regard the numerical solution in $H_{64}$ as a reference 'real' solution to detect weak errors and corresponding convergence orders for $\gamma=0$ and $\gamma=\frac{1}{8}$. The errors versus $N$ are plotted in log-log scale in Figure \ref{figure14.1}. We observe that the average weak convergence order is consistent with predicted theoretical rate.

\bibliographystyle{plain}
\bibliography{references.bib}

\end{document}